\documentclass[11pt]{article}

\usepackage{amsmath,amssymb,amsthm,amscd,mathrsfs}
\usepackage{latexsym,amsmath,amssymb,graphicx}
\usepackage[all]{xy} 
\usepackage{multirow}
\usepackage{color}
\usepackage{epsf}
\usepackage{youngtab}
\xyoption{v2} \xyoption{2cell} \xyoption{dvips}
\CompileMatrices \numberwithin{equation}{section}
\newtheorem{prop}{Proposition}[section]
\newtheorem{theo}[prop]{Theorem}
\newtheorem{lemm}[prop]{Lemma}
\newtheorem{coro}[prop]{Corollary}
\newtheorem{rema}[prop]{Remark}

\newtheorem{defi}[prop]{Definition}

\numberwithin{equation}{section}

\newcommand{\be}{\begin{equation}}
\newcommand{\ee}{\end{equation}}

\newcommand{\IP}{\mathbb{P}}

\newcommand\IZ{\mathbb {Z}}

\newcommand\IQ{\mathbb {Q}}
\newcommand{\IC}{\mathbb{C}}

\newcommand{\IR}{\mathbb{R}}

\newcommand{\ba}{\begin{array}}
\newcommand{\ea}{\end{array}}

\newcommand{\CV}{{\mathcal V}}

\newcommand{\IF}{{\mathbb F}}

\newcommand{\wG}{{\widetilde G}}

\newcommand{\CS}{{\mathcal S}}

\newcommand{\CK}{{\mathcal K}}
\newcommand{\CN}{{\mathcal N}}

\newcommand{\IH}{{\mathbb H}}
\newcommand{\bal}{\begin{aligned}}
\newcommand{\eal}{\end{aligned}}

\newcommand{\hX}{{\widehat X}}
\newcommand{\hZ}{{\widehat Z}}

\newcommand{\CZ}{{\mathcal Z}}

\newcommand{\longto}{\longrightarrow}

\newcommand{\ch}{{\mathrm{ch}}}

\newcommand{\CO}{{\mathcal O}}

\newcommand{\CH}{{\mathcal H}}
\newcommand{\CF}{{\mathcal F}}

\newcommand{\CM}{{\mathcal M}}

\newcommand{\CI}{{\mathcal I}}

\newcommand{\CQ}{{\mathcal Q}}

\newcommand{\wS}{{\widetilde S}}

\newcommand{\CT}{{\mathcal T}}

\newcommand{\calP}{{\mathcal P}}
\newcommand{\wXi}{{\widetilde \Xi}}
\newcommand{\tbeta}{{\tilde \beta}}

\newcommand{\wF}{{\widetilde F}}

\newcommand{\hF}{{\widehat F}}
\newcommand{\hG}{{\widehat G}}

\newcommand{\hE}{{\widehat E}}

\newcommand{\hPhi}{{\widehat \Phi}}
\newcommand{\hU}{{\widehat U}}
\CompileMatrices \LaTeXdiagrams \UseAllTwocells
\newdimen\tableauside\tableauside=1.0ex
\newdimen\tableaurule\tableaurule=0.4pt
\newdimen\tableaustep
\def\phantomhrule#1{\hbox{\vbox to0pt{\hrule height\tableaurule width#1\vss}}}
\def\phantomvrule#1{\vbox{\hbox to0pt{\vrule width\tableaurule height#1\hss}}}
\def\sqr{\vbox{%
  \phantomhrule\tableaustep
  \hbox{\phantomvrule\tableaustep\kern\tableaustep\phantomvrule\tableaustep}%
  \hbox{\vbox{\phantomhrule\tableauside}\kern-\tableaurule}}}
\def\squares#1{\hbox{\count0=#1\noindent\loop\sqr
  \advance\count0 by-1 \ifnum\count0>0\repeat}}
\def\tableau#1{\vcenter{\offinterlineskip
  \tableaustep=\tableauside\advance\tableaustep by-\tableaurule
  \kern\normallineskip\hbox
    {\kern\normallineskip\vbox
      {\gettableau#1 0 }%
     \kern\normallineskip\kern\tableaurule}%
  \kern\normallineskip\kern\tableaurule}}
\def\gettableau#1 {\ifnum#1=0\let\next=\null\else
  \squares{#1}\let\next=\gettableau\fi\next}
\begin{document}

\title{Vertical sheaves and Fourier-Mukai transform on elliptic 
Calabi-Yau threefolds}
\author{D.E. Diaconescu}
\date{}
\maketitle

\begin{abstract} 
This paper studies the action of the Fourier-Mukai transform on moduli 
spaces of vertical 
torsion  sheaves on elliptic Calabi-Yau threefolds in Weierstrass form. 
Moduli stacks of semistable one dimensional sheaves on such threefolds are identified with open and closed substacks
of moduli stacks of vertical semistable  two dimensional sheaves on 
their Fourier-Mukai duals. In particular, this yields explicit conjectural results for  
Donaldson-Thomas invariants of vertical two dimensional sheaves on 
K3-fibered elliptic Calabi-Yau threefolds. 
\end{abstract}

\tableofcontents

\section{Introduction}\label{intro}

Starting with Mukai's work on the subject \cite{FMfunctor, semihombundles}, Fourier-Mukai functors 
have played a central role in the study of moduli spaces of stable sheaves
on algebraic varieties. An incomplete list of applications of Fourier-Mukai 
functors to 
moduli spaces of torsion free sheaves on surfaces includes 
\cite{FM_stab_K3, FM_mirror, FM_existence, FM_elliptic, Mod_sheaves_abelian,FM_transf_abelian,FM_instantons,St_sheaves_elliptic,
Semihom_abelian,
Stab_FM_transform_I,Stab_FM_transform_II}.
Further applications to moduli spaces of torsion free sheaves on 
elliptic threefolds and higer dimensional elliptic fibrations include \cite{Vb_F_theory,Vb_elliptic,FM_K3_elliptic,Relstabelliptic,Fiberwise,
Twisted_elliptic,St_sheaves_K3_fibr}.
A comprehensive review of the subject and a more complete list of results 
can be be found in \cite{FM_book}. 
More recently, t-structures and moduli problems of Bridgeland stable objects in the derived category 
have been studied in 
\cite{FM_wallcrossing,mod_complex_Kthree,Some_moduli,Rank_one_stable, FM_perv_coh,EFM_algorithm,Stab_FM_elliptic, StabFM_II,Elliptic_t_structures} using a similar approach.  

Of particular importance 
for the present paper is the relative 
Fourier-Mukai transform for elliptic fibrations. This was constructed by 
Bartocci et al \cite{FM_mirror} and Bridgeland \cite{FM_elliptic, FM_thesis} for elliptic surfaces and 
Friedman, Morgan and Witten  \cite{Vb_F_theory,Vb_elliptic} 
for stable bundles on elliptic threefolds. 
The foundational results for elliptic threefolds used in this paper were proven by Bridgeland and Macciocia in \cite{FM_K3_elliptic}. The higher 
dimensional construction was carried out by Bartocci et al in 
\cite{Relstabelliptic}. 

An important problem in this framework is whether the Fourier-Mukai 
transform preserves Gieseker stability,  in particular if it yields 
isomorphisms of moduli spaces of semistable sheaves. Several results 
obtained in the literature prove that this is the indeed the case for suitable open subspaces of moduli 
spaces parameterizing relatively semistable objects. However  isomorphisms 
of proper moduli spaces are much harder to prove.  One such 
result was obtained by Yoshioka in \cite{Mod_sheaves_abelian}, showing that Fourier-Mukai transform identifies moduli spaces 
of semistable pure dimension one sheaves on an elliptic surface 
with moduli spaces of semistable torsion-free sheaves on the dual 
surface. The main goal of the present paper is to study the analogous 
problem for pure dimension one sheaves on elliptic threefolds. As explained 
in more detail below this problem is mainly motivated by applications 
to Donaldson-Thomas invariants \cite{RT,genDTI} of pure dimension two sheaves and 
modularity questions.

\subsection{The main result}\label{mainres} 
Let $p:X\to B$ be a smooth projective Weierstrass model with trivial 
canonical class 
over a smooth Fano surface $B$. The Mukai dual $\hX$ of $X$ was 
constructed in \cite{FM_K3_elliptic} as a fine relative moduli space 
for rank one degree zero torsion free sheaves on the fibers of $p:X\to B$. 
For sufficiently generic $X$ the dual $\hX$ is again a smooth Weietrstrass 
model ${\hat p}: \hX \to B$ and there is a canonical isomorphism 
$\hX \simeq X$ over $B$. Since $\hX$ is a fine moduli space, there is a 
(non-unique) 
universal Poincar\'e sheaf $\cal P$ on $\hX\times X$. The Fourier-Mukai 
functor $\Phi : D^b(\hX)\to D^b(X)$ with kernel $\calP$ was proven in \cite{FM_K3_elliptic} to be an 
equivalence of derived categories. Moreover it was also shown there in that the inverse functor  
$\hPhi: D^b(X)\to D^b(\hX)$ is also a Fourier-Mukai transform whose 
kernel $\CQ$ is the derived dual $\calP$ up to a shift.  A more detailed summary is provided in Section \ref{FMbasics}.

The main goal of this paper is to study the action of the above 
Fourier-Mukai 
functors on moduli stacks of Gieseker semistable torsion sheaves on $X,\hX$. 
In order to formulate a concrete statement one first needs a concrete presentation of the K\"ahler cones and the homology groups of $X, \hX$. As shown in Lemma \ref{basicX}, 
one has an isomorphism \[
{\rm Pic}(X)/_{\rm torsion} \simeq \IZ\langle \Theta \rangle  \oplus p^*{\rm Pic}(B),
\] 
where $\Theta$ is the image of the canonical section $\sigma: B \to X$. 
Then any K\"ahler class $\omega\in {\rm Pic}_\IR(X)$ can be written as $\omega = 
t \Theta + p^*\eta$, with $t\in \IR$, $t>0$ 
and $\eta \in {\rm Pic}_\IR(B)$ a sufficiently 
ample K\"ahler class on $B$. In particular $\omega = t \Theta -s p^*K_B$ 
is a K\"ahler class on $X$ for $s>t>0$, where 
$K_B$ is the canonical class of $B$. Lemma \ref{basicX} also shows that there is a natural isomorphism ${\rm Pic}(X) \simeq H_4(X,\IZ)$ which will be used implicitely throughout this paper. In particular the pairing between 
K\"ahler classes and homology classes will be identified with the intersection product. Using Poincar\'e duality, Chern classes of sheaves on $X$ will be also regarded as even homology classes. 
Finally, note the  direct sum decomposition 
\[
H_2(X,\IZ)/_{\rm torsion} \simeq \IZ\langle f \rangle \oplus \sigma_*H_2(B,\IZ) 
\]
where $\sigma:B\to X$ is the canonical section of the Weierstrass 
model and $f$ is the elliptic fiber class. 
Of course, completely analogous statements 
hold for ${\hat p}: \hX \to B$, the notation being obvious. 

This paper will concrentrate on the relation between pure dimension 
one sheaves on $\hX$ and vertical pure dimension two sheaves on $X$. According to Definition \ref{vertdef}.i, a sheaf $E$ on $X$ of 
pure dimension two is vertical if $\ch_1(E)\cdot f=0$ and 
$\ch_2(E)$ is a multiple of the fiber class $f$.  
The discrete invariants of a sheaf $\hF$ of pure dimension one on $\hX$ are given by an element 
\[
{\hat \gamma} = ({\hat \gamma}_i)_{1\leq i \leq 3} 
\in H_2(B,\IZ)\oplus \IZ\oplus \IZ
\]
where 
\[ 
\ch_2(\hF)= {\hat \sigma}_*({\hat \gamma}_1)+ 
{\hat \gamma}_2 {\hat f}, \qquad 
\chi(\hF) = {\hat \gamma}_3.
\]
The discrete invariants of a vertical  sheaf $E$ on $X$ 
of pure dimension two are given by 
\[ 
\gamma=(\gamma_i)_{1\leq i\leq 3}\in H_2(B,\IZ) \oplus 
(1/2)\IZ \oplus \IZ,
\] 
where
\[
 \ch_1(E)=p^*\gamma_1, \qquad \ch_2(E)=\gamma_2 f, \qquad 
\ch_3(E) = -\gamma_3 \ch_3(\CO_x) 
\]
with $x\in X$ an arbitrary closed point. 
According to  equations \eqref{eq:FMchernA}, 
\eqref{eq:FMchernB}, the induced action of Fourier-Mukai transform on 
numerical invariants is encoded in the 
group isomorphism 
\[
\bal 
\phi: H_2(B,\IZ)\oplus \IZ  \oplus & \IZ\, {\buildrel \sim \over \longto}\,  
H_2(B,\IZ) \oplus (1/2)\IZ \oplus \IZ, \\
 \phi({\hat \gamma}_1,{\hat \gamma}_2, {\hat \gamma}_3) & = (
{\hat \gamma}_1, {\hat \gamma}_3 + K_B\cdot {\hat \gamma}_1/2, 
{\hat \gamma}_2).\\
\eal
\]
Here $\cdot$ denotes the intersection product on $B$. 

Note also that Definition \ref{vertdef}.ii
introduces a notion of adiabatic stability for vertical  sheaves on $X$ which plays an important part in this paper. Given a K\"ahler class $\omega = t\Theta + p^*\eta$, a vertical pure dimension two 
sheaf $E$ is called $\omega$-adiabatically semistable if and only if it is  Gieseker semistable with respect to all K\"ahler classes 
$\omega' = t' \Theta + p^*\eta$, where $0< t'\leq t$. 

Given K\"ahler classes $\omega = t\Theta - sp^*K_B$, ${\hat \omega} = {\widehat \Theta} - s{\hat p}^*K_B$ with $s>t>0$, $s>1$, 
let $\CM_{\hat \omega}(\hX, {\hat \gamma})$, $\CM_{\omega}(X,\gamma)$ denote the moduli stacks
of Gieseker semistable sheaves with discrete invariants ${\hat \gamma}$, $\gamma$ on $\hX$, $X$, respectively. Let $\CM^{\sf ad}_{\omega}(X,\gamma)\subset \CM_{\omega}(X,\gamma)$ 
be the substack of adiabatically semistable sheaves as defined 
in \ref{vertdef}.ii. 
Then the main result of the present paper is 
\begin{theo}\label{main} 
Let 
${\hat \gamma} \in H_2(B,\IZ)\oplus \IZ\oplus \IZ$ be fixed numerical 
invariants such that ${\hat \gamma}_3 >0$.
Then there exists a constant $s_1({\hat \gamma})\in \IR$, 
$s_1({\hat \gamma})>1$, depending 
on ${\hat \gamma}$, such that 
for any $s\in \IR$, $s> s_1({\hat \gamma})$, there exists a second constant $t_1({\hat \gamma},s)\in \IR$, 
$0<t_1({\hat \gamma},s)<1$, depending 
on $({\hat \gamma},s)$, such that the following  statements hold 
for any $t\in \IR$, $0<t<t_1(s, {\hat \gamma})$. 

$(i)$ The Fourier-Mukai transform $\Phi$ yields an isomorphism 
of moduli stacks 
\[ 
\varphi: \CM_{\hat \omega}(\hX, {\hat \gamma})
{\buildrel \sim \over \longto} \CM^{\sf ad}_{\omega}(X,\gamma), 
\]
where ${\hat \omega} = \Theta - s{\hat p}^*K_B$, 
$\omega = t\Theta -s p^*K_B$ and $\gamma=\phi({\hat \gamma})$. 

$(ii)$ The substack 
$\CM^{\sf ad}_{\omega}(X,\gamma)\subset \CM_{\omega}(X,\gamma)$ is open and closed in $\CM^{}_{\omega}(X,\gamma)$.
\end{theo}

The proof of Theorem \ref{main} is given in Section \ref{FMvertical} and 
requires some preliminary results proven in Section \ref{verticalsect}. 
In comparison with the analogous result for elliptic surfaces 
\cite[Thm. 3.15]{Mod_sheaves_abelian}, one needs to introduce a suitable notion 
of generic stability  for vertical pure dimension two sheaves in Definition  \ref{genstabdef}. Then one has to further check that generic stability 
is equivalent to adiabatic stability in Lemmas \ref{genstablemma}
and \ref{adiabstablemma}. 
The proof is then given step-by-step in Section \ref{FMvertical}. 
In contrast with \cite[Thm 3.15]{Mod_sheaves_abelian}, one cannot rule 
out non-adiabatic components of the moduli stack of semistable pure dimension two sheaves on a threefold 
by taking an appropriate limit of the K\"ahler class. 
However, as shown below, such components can be ruled out 
for elliptic threefolds which also admit a $K3$-fibration 
structure, and for two dimensional sheaves supported 
on the $K3$ fibers.

\subsection{Sheaf counting on elliptic $K3$ pencils}\label{countingsect}

As stated in the second paragraph of the introduction, Theorem 
\ref{main} is mainly motivated by applications to Donaldson-Thomas 
invariants of pure dimension two sheaves on elliptic Calabi-Yau threefolds. 
These are counting invariants defined in \cite{RT} for stable sheaves 
and generalized in \cite{genDTI,wallcrossing} for semistable ones. 
Generating 
series of  Donaldson-Thomas invariants for pure dimension two sheaves 
have been conjectured to have modular properties in \cite{Denef:2007vg, 
M5_genus}. In the mathematics literature, this conjecture 
has been proven for certain cases in 
\cite{DT_twodim,Gen_DT_twodim,genDT_local_plane}. In particular 
explicit results for Donaldson-Thomas invariants of such sheaves on $K3$ fibered Calabi-Yau threefolds were obtained by Gholampour and Sheshmani 
in \cite{DT_twodim}. 
For nodal $K3$ pencils these results are restricted to rank one 
torsion free sheaves on reduced $K3$ fibers. 

On the other hand, string theoretic arguments \cite{E_strings,Quantum_elliptic} lead to a conjectural identification 
of Donaldson-Thomas invariants for vertical pure dimension two 
sheaves on an elliptic threefold $X$ with 
genus zero Gopakumar-Vafa invariants on its dual $\hX$. 
This correspondence was first conjectured in \cite{E_strings} for 
sheaves supported on a rational elliptic surface inside $X$. 
As observed in \cite{HST}, in that case this follows from the 
results of \cite{Mod_sheaves_abelian}. 
As it stands, Theorem \ref{main} does not prove such an 
identification for general vertical sheaves 
because the moduli stack $\CM_\omega(X,\gamma)$ 
can in principle have other components in addition 
to $\CM_{\omega}^{\sf ad}(X,\gamma)$. From a string theory 
point of view it is natural to conjecture that such components are absent for sufficiently small $t_1({\hat \gamma},s)$, but mathematically this is an open problem. 

As shown in below, there is however one situation where such extra components can be ruled out. 
Excepting $\IP^2$, all smooth Fano surfaces $B$ have a natural 
projection $\rho: B \to \IP^1$, which induces a projection $\pi=\rho\circ p : X \to \IP^1$. The generic fiber of $\rho$ is a 
smooth reduced elliptic $K3$ surface on $X$. Moreover if $B$ is a 
Hirzebruch surface $\IF_a$, $0\leq a\leq 1$, for sufficiently 
generic $X$, all fibers are reduced irreducible $K3$ surfaces with at most nodal singularities. Under this assumptions, 
Proposition \ref{Kthreepencil} shows that no extra components 
are present in the moduli space of semistable vertical sheaves supported on $K3$ fibers for suitable K\"ahler classes. 
Therefore, in such cases Theorem \ref{main} 
yields explicit conjectural results for generalized  
Donaldson-Thomas invariants of two dimensional 
sheaves supported on the $K3$ fibers, 
verifying the modularity conjecture. 

In more detail, suppose $B$ is the total space of the projective bundle  $\IP(\CO_{\IP^1} 
\oplus \CO_{\IP^1}(a))$ with $0\leq a \leq 1$. Let $\rho: B \to \IP^1$ and $\pi = \rho \circ p : X \to \IP^1$. For sufficiently generic $X$ the fibers of 
$\pi$ are reduced irreducible $K3$-surfaces with at most nodal singularities. 
Recall that  ${\hat \sigma} : B \to \hX$ denotes the canonical section
and 
${\hat f}$ denotes the fiber class of ${\hat p}: \hX \to B$. Let 
$\Xi$ denote the fiber class of $\rho:B \to \IP^1$. 
The $K3$ fiber class on $X$ is $D=p^*\Xi$. 
Then note the following. 

\begin{prop}\label{Kthreepencil} 
Let $\omega = 
t\Theta - sp^*K_B$ with  $s,t \in \IR$, $s>t >0$.  Let 
$ \gamma = (rD, l, m)$ be arbitrary discrete invariants with 
$r,l,m\in \IZ$, $r\geq 1$. 
Then there exists a constant $t_2(\gamma, s)\in \IR$, 
$t_2(\gamma, s)>0$ such that $\CM^{\sf ad}_{\omega}(X,\gamma)= \CM_{\omega}(X,\gamma)$ for any $0<t< t_2(\gamma,s)$ such that  
$t/s\in \IR\setminus \IQ$.   
\end{prop}

The proof of Proposition \ref{Kthreepencil} is given in Section 
\ref{verticalKthree}. It should be noted that similar results for 
torsion free 
sheaves on elliptic surfaces have been obtained before in 
\cite[Thm. I.3.3]{vb_inv_elliptic}, \cite[Prop. I.1.6]{weight_two_hodge} and \cite[Lemma 1.2]{notes_elliptic}. Here one has to generalize 
these results to semistable pure dimension two sheaves supported 
on scheme theoretic thickenings of divisors in the $K3$ pencil, including nodal fibers. This requires a careful reduction to the reduced smooth surface 
case via Jordan-H\"older filtrations and blow-ups.

Next consider numerical invariants ${\hat \gamma} = (r\Xi,n , k)$
in Theorem \ref{main},
where $r,n,k\in \IZ$, $r,k\geq 1$, $n\geq 0$.
Then equations \eqref{eq:FMchernB} yield 
\be\label{eq:twodinvariants}
\gamma_1 = rD, \qquad \gamma_2 = (k-r)f, \qquad \gamma_3 = -n. 
\ee
Let $DT_{\hat \omega}(\hX; r,n,k)\in \IQ$ denote the generalized Donaldson-Thomas 
invariants counting ${\hat \omega}$-semistable pure dimension one sheaves on $\hX$ constructed in \cite{genDTI}. According to \cite[Thm 6.16.a]{genDTI}, 
these invariants are independent of ${\hat \omega}$, hence the subscript 
will be dropped in the following. Moreover it is conjectured in 
 \cite[Conj. 6.12]{genDTI} that there exist integral invariants 
$\Omega(\hX; r,n,k)\in \IZ$ related to the rational ones by the multicover formula 
\be\label{eq:multicoverA}
DT(\hX; r,n,k)=\sum_{\substack{m\in \IZ,\ m\geq 1,\\ m|(r,n,k)}}
{1\over m^2}\Omega(\hX; r/m,n/m,k/m).
\ee
Alternatively the integral invariants can be conjecturally 
defined directly 
by specialization of the motivic invariants of Kontsevich and Soibelman \cite{wallcrossing}
as explained in Section 7.1 of loc.cit. 

For a primitive vector ${\hat \gamma}= 
(r,n,k)$ there are no strictly semistable objects, and the invariants 
$DT(\hX; r,n,k)=\Omega(\hX; r,n,k)$ specialize to the integral virtual cycle invariants 
defined in \cite{RT}. In particular this holds for $k=1$. Then the 
resulting invariants were conjecturally identified with genus zero Gopakumar-Vafa invariants in \cite{genuszero}, 
\be\label{eq:GVDT}
\Omega(\hX, r,n,1) = n_0(\hX, r,n)
\ee
for any $(r,n)\in \IZ^2$, $r,n\geq 0$, $(r,n)\neq (0,0)$. 
Here $n_0(r,n)$ denotes the genus zero Gopakumar-Vafa invariant for curve 
class $r{\hat\sigma}(\Xi) + n{\hat f}$. 
As shown in \cite{stab_counting}, equation \eqref{eq:GVDT} 
follows from the GW/stable pair correspondence conjectured 
in \cite{stabpairsI} provided that the integral 
invariants 
$\Omega(\hX; r,n,k)\in \IZ$ are independent of $k$
for fixed $(r,n)$. Independence of $k$ is a special 
case of \cite[Conj. 6.13]{stab_counting}. 

Let $DT_\omega(X,r,l,m)$ denote the generalized Donaldson-Thomas invariants counting $\omega$-semistable vertical 
two dimensional sheaves 
with invariants $\gamma=(r\Xi,l,m)$ . The wallcrossing formulas 
of \cite{genDTI, wallcrossing} imply easily that the invariants $DT_\omega(X,r,l,m)$ are independent of $\omega$, hence the 
subscript may be dropped. Again, Conjecture 6.12 in  \cite{genDTI} 
states the existence of integral invariants 
$\Omega(X,r,l,m)$ 
related to the rational ones by a multicover formula of the form
\eqref{eq:multicoverA}.

Then Theorem \ref{main} and Proposition 
\ref{Kthreepencil} imply that 
\[
DT(\hX; r,n,k) = DT(X, r, k-r,n) 
\]
for any $r,n,k\in \IZ$, $r,k\geq 1$, $n\geq 0$. 
Granting the existence of integral invariants, they will be also related by 
\[ 
{\Omega}(\hX; r,n,k) = {\Omega}(X, r, k-r,n).
\]
However note that there is an isomorphism of moduli stacks 
$\CM_\omega(X,r,l,m) \simeq \CM_\omega(X,r,l-r,m)$ for any 
$(r,l,m)$. This is obtained by taking tensor product by the line bundle 
$p^*\CO_B(-C_0)$, where $C_0$ is a section of the ruling $\rho:B\to \IP^1$. 
For concreteness let $C_0$ be the unique section with $C_0^2=-1$ for 
$B = \IF_1$ and an arbitrary section with $C_0^2=0$ 
for $B = \IF_0$. 
Therefore 
\be\label{eq:FMDT}
{\Omega}(\hX; r,n,k) = {\Omega}(X, r, k,n).
\ee
for any $(r,n,k)$, $r,k\geq 1$, $n\geq 0$. 

Now let 
\[
Z_{X,r,k}(q) = \sum_{n\in \IZ}  {\Omega}(X, r, k,n) q^{n-r/2}. 
\]
and suppose for concreteness that $B=\IF_1$. Then, granting 
the invariance of $DT(\hX, r,n,k)$ under translations $k\mapsto k+1$ 
and the identification \eqref{eq:GVDT} one obtains 
\[ 
Z_{X,r,k}(q) = \sum_{n\geq 0} n_0(\hX, r,n) q^{n-r/2} 
\]
for any $r,k\in \IZ$, $r,k\geq 1$.

Granting the identification \eqref{eq:GVDT}, an explicit formula for 
the series $Z_{X,r,k}(q)$ follows from the work of Maulik and Pandharipande \cite{GW_NL} on Gopakumar-Vafa invariants 
of K3 pencils. The explicit computation for the Weierstrass model
over $\IF_1$ was done  
by Rose and Yui \cite{GV_elliptic}. 
The formula obtained in \cite[Thm. 7.5]{GV_elliptic} is written in 
terms of a certain transformation of modular forms 
defined in 
\cite[Def. 7.1]{GV_elliptic}.  Let  
\[ 
f(z) = \sum_{n} a_n z^n 
\]
 be modular form for $SL(2,\IZ)$ and $r,k\in \IZ$, $r\geq 1$. 
Then define $f_{r,k}(z)$ by 
\[
f_{r,k}(z) = \sum_{n=0}^\infty a_{rn+k'} z^{rn+k'}
\]
where $0\leq k'< r$ is the unique integer in this range such 
that $k' \equiv k$ (mod $r$). Note
that this is modular form 
for the subgroup $\Gamma_1(r^2) 
\subset SL(2,\IZ)$ of the same weight as $f(z)$.
Then identity \eqref{eq:FMDT} and \cite[Thm 7.5]{GV_elliptic} yield
the following conjectural formula:
\be\label{eq:Zconj}
Z_{X,r,k}(q) = -2 \sum_{\ell=0}^{r-1} \left({1\over \Delta(u)}\right)_{r,\ell-1} E_{10}(u)_{r,1-\ell}
\ee
where $q=u^r$ and $\Delta(u) = {1/ \eta(u)^{24}}$. 

To conclude, note two natural open problems emerging from the present work. 
One open question in the  context of Theorem \ref{main} is whether there 
exists a sufficiently small constant $t_1({\hat \gamma}, s)$ such that the 
moduli stack 
$\CM_\omega(X,\gamma)$ coincides with the 
substack of adiabatically semistable objects. 
String theoretic 
arguments \cite{E_strings, Quantum_elliptic} lead to the conjecture 
that this is indeed the case. 
As shown in Proposition 
\ref{Kthreepencil}, this holds in the special case of vertical sheaves 
on elliptic K3 pencils. The proof given in Section \ref{verticalKthree} 
relies on Bogomolov inequality and the algebraic Hodge theorem for 
surfaces. This leads to the interesting question whether analogous tools can be developed in general for vertical sheaves on elliptic 
threefolds. 

The second open problem is whether formula \eqref{eq:Zconj} 
can be given a direct proof using degeneration techniques as in \cite{DT_twodim}. 

{\it Acknowledgements.} This work is motivated by a related  project initiated in collaboration with Vincent Bouchard, 
Thomas Creutzig, Chuck Doran, Terry Gannon and Callum Quigley. 
I am very grateful to them for very stimulating conversations and sharing their insights with me. I would also like to thank Lev Borisov, Ugo Bruzzo, Ron Donagi, 
Tony Pantev and Noriko Yui for very helpful discussions. The author 
also acknowledges the partial support of  NSF grant DMS-1501612 during the completion if this work.

\section{Vertical sheaves and adiabatic stability}\label{verticalsect}

This section introduces adiatically semistable vertical sheaves on 
elliptic Weierstrass models, and shows that adiabatic stability is equivalent to a natural notion of generic stability.

\subsection{Basics of Weirestrass models}
Let $B$ be a smooth projective del Pezzo surface. 
Let $p:X\to B$ be a smooth generic Weierstrass model  with canonical section $\sigma: B \to X$.  
Let $\Theta\subset X$ denote the image of the canonical section. 
Then $\Theta$ determines a homology class in $H_4(X,\IZ)$ as well as a 
divisor class in ${\rm Pic}(X)$. 
Let $f\in H_2(X,\IZ)$ denote the class of the elliptic fiber. 
The same notation $\cdot$ will be used for 
the intersection product on $X$, as well as $B$. The distinction will be clear from the context. 

\begin{lemm}\label{basicX}
There are direct sum decompositions 
\be\label{eq:Xhomology} 
H_4(X,\IZ)/_{\rm torsion} \simeq \IZ\langle \Theta\rangle \oplus p^*H_2(B,\IZ)\qquad 
H_2(X,\IZ)/_{\rm torsion} \simeq \IZ\langle f\rangle \oplus \sigma_*H_2(B,\IZ), 
\ee
Moreover, there is an isomorphism ${\rm Pic}(X)\simeq H^4(X,\IZ)$.
\end{lemm} 

{\it Proof.} One proceeds by analogy with 
\cite[Lemma 6.1]{GV_elliptic}.
Note that $h^{0,i}(X)=0$, $i\in \{1,2\}$, and $h^{1,1}(X) = h^{1,1}(B)+1$ 
according to 
\cite[Sect. 11]{MW_rank}. This implies that there is an isomorphism 
${\rm Pic}(X) \simeq H^2(X,\IZ)$. By Alexander-Lefschetz duality, 
there is also an isomorphism $H^2(X,\IZ) \simeq H_4(X,\IZ)$. 
Next recall that 
${\rm Pic}(B)\simeq H_2(B,\IZ)$ is freely generated 
by rational curve classes $C_1, \ldots, C_{h^{1,1}(B)}$ such that the intersection 
matrix $I_B=(C_i\cdot C_j)_{1\leq i,j \leq h^{1,1}(B)}$ has determinant 
$|{\rm det}(I_B)|=1$. Let $D_i = p^*C_i\in {\rm Pic}(X)\simeq H_4(X,\IZ)$, $1 \leq i \leq h^{1,1}(B)$.
Let $I_X$ denote the intersection matrix 
between the divisor classes $\Theta, D_1, \ldots, D_{h^{1,1}(B)}$ and 
the curve classes $f, \sigma_*(C_1), \ldots, \sigma_*(C_{h^{1,1}(B)})$
on $X$. Straightfowrard intersection computations show that $|{\rm det}(I_X)|=1$ as well. 
This implies the isomorphisms claimed above.

\hfill $\Box$

As explained in Section \ref{mainres}, the isomorphism ${\rm Pic}(X)\simeq H_4(X,\IZ)$ following from Lemma \ref{basicX} will be 
implicitely used throughout this paper. Moreover, Chern 
classes of sheaves on $X$ will be identified with homology classes 
by Poincar\'e duality.  Then note the following.
\begin{coro}\label{basicXcoro} 
$(i)$ A real divisor class 
\[ 
\omega = t \Theta + p^*\eta, \qquad t\in \IR, \ t>0, 
\]
is ample if and only if $\eta + tK_B$ is an ample 
divisor class on $B$. 

$(ii)$ Let $C\in H_2(B,\IZ)$ be an arbitrary 
curve class and let $\Sigma = \sigma_*(C) + nf \in H_2(X,\IZ)/_{\sf torsion}$ with $n \in \IZ$. Then $\Sigma$ is an efective curve class 
if and only if $C$ is effective and $n\geq 0$.
\end{coro} 

{\it Proof.} For $(i)$ suppose $\Sigma$ is an effective curve class 
on $X$ which contains an irreducible curve. 
Let $\eta\in {\rm Pic}(B)$ be an ample class and note that 
\[
\Sigma \cdot p^*\eta = p_*\Sigma \cdot \eta \geq 0. 
\]
Since $\Sigma$ contains an irreducible curve, one of the following cases must hold.
\begin{itemize}
\item[$(a)$] The set theoretic support of the irreducible curve in 
$\Sigma$ is not contained in $\Theta$. In this case 
$\Sigma\cdot \Theta\geq 0$. 
\item[$(b)$] The set theoretic support of the irreducible curve in 
$\Sigma$ is contained in $\Theta$. In this case
 $\Sigma \cdot \Theta <0$ and $\Sigma =\sigma_*(C)$ with $C$ an effective 
curve class on $B$. Moreover, 
\[ 
\Sigma \cdot \Theta = C\cdot K_B.
\]
\end{itemize}
Then the claim $(i)$ follows easily. 

$(ii)$ Let $\eta \in {\rm Pic}(B)$ be an arbitrary ample class. Note that 
for sufficiently large $k>0$ there exists a divisor $H$ in the linear system 
$|k\eta|$ such that  $Z=p^{-1}(C)$ does not contain the set theoretic 
support of any of the irreducible components of $\Sigma$. 
Since $\Sigma$ is effective, this implies $\Sigma \cdot \eta \geq 0$, 
hence 
\[
C\cdot \eta =  p_*\Sigma \cdot \eta=\Sigma \cdot p^*\eta  \geq 0. 
\]
Since this holds for any ample class $\eta$, it follows that $C$ must be effective or zero. If $C=0$, the claim is obvious. Suppose $C\neq 0$ and  
$n<0$. Note that $K_B\cdot C\neq 0$. 
Let $s\in \IR$ be a real number such that 
\be\label{eq:snumber}
0< s-1 <  {|n|\over |K_B \cdot C|}. 
\ee
Then $\eta = \Theta - sp^*K_B$ is an ample class on $X$, hence 
\[
0< \eta\cdot \Sigma = (s-1) |K_B\cdot C| -|n|. 
\]
This contradicts the second inequality in \eqref{eq:snumber}.

\hfill $\Box$.

\subsection{Adiabatic and generic stability}\label{adgensect}

Recall Gieseker and slope 
stability for two dimensional sheaves on $X$.
Let $\omega$ be an ample 
class on $X$.  
For any nonzero coherent sheaf $E$ on 
$X$ of dimension two let 
\[
\mu_{\omega}(E) = {\omega \cdot \ch_2(E)\over 
\omega^2\cdot \ch_1(E)/2}, \qquad
\nu_{\omega}(E) = {\chi(E)\over \omega^2\cdot  \ch_1(E) /2}.
\]
Then Gieseker (semi)stability  with respect to
$\omega$ is defined by
the conditions
\be\label{eq:twistedstabA}
\mu_{\omega}(E') \ (\leq) \ \mu_{\omega}(E)
\ee
for any proper nonzero
subsheaf $0\subset E' \subset E$, and
\be\label{eq:twistedstabB}
\nu_{\omega}(E') \ (\leq) \ \nu_{\omega}(E)
\ee
if the slope inequality \eqref{eq:twistedstabA} is saturated. 
Recall that any Gieseker semistable sheaf must be of pure dimension. 
Furthermore, a pure dimension two sheaf $E$ is Gieseker semistable 
if and only if the above inequalities are satisfied for saturated proper 
nonzero subsheaves i.e. $E/E'$ pure of dimension two. 
In contrast, $\omega$-slope (semi)stability is defined by imposing only 
condition \eqref{eq:twistedstabA} with respect to nonzero proper 
saturated subsheaves. 

For completeness, recall that the $\omega$-slope of
a nonzero 
pure dimension one sheaf $E$ is defined by 
\[
\mu_\omega(E) = {\chi(E)\over \omega\cdot \ch_2(E)}.
\]
Such a sheaf is 
called 
Gieseker $\omega$-semistable if and only if 
\[
\mu_\omega(E') \ (\leq)\ \mu_\omega(E)
\]
for any proper nontrivial subsheaf $E'\subset E$. 
In this case Gieseker $\omega$-(semi)stability and $\omega$-slope (semi)stability coincide. 

Throughout 
this paper Gieseker stability relative to an ample class $\omega\in {\rm Pic}_\IR(X)$ will be simply called  $\omega$-stability.

\begin{defi}\label{vertdef}
$(i)$ 
A pure dimension two sheaf $E$ on $X$ will be called vertical 
if and only if 
\[ 
\ch_1(E) \in p^*{\rm Pic}(B), \qquad \ch_2(E) \in (1/2)\IZ\langle f\rangle. 
\] 

$(ii)$ A vertical pure dimension two sheaf $E$ on $X$ will be called 
adiabatically $\omega$-(semi)stable 
if and only if 
it is $(t'\Theta + p^*\eta)$-(semi)stable for all 
$0 < t' \leq t$. 

$(iii)$ A vertical pure dimension two sheaf $E$ on $X$ will be called 
adiabatically $\omega$-slope (semi)stable 
if and only if 
it is $(t'\Theta + p^*\eta)$-slope (semi)stable for all 
$0 < t' \leq t$. 

\end{defi} 

Note that the discrete invariants of a vertical sheaf $E$ are given by 
\be\label{eq:numinvA}
\ch_1(E)=p^*C, \qquad \ch_2(E)=kf, \qquad \ch_3(E)= -n \ch_3(\CO_x) 
\ee
where $C\in {\rm Pic}(B)$ is an effective divisor class on $B$, 
$k\in  (1/2)\IZ$ and $n\in \IZ$. Using the isomorphism ${\rm Pic}(B) 
\simeq H_2(B,\IZ)$, this yields an element 
\[ 
\gamma = (C,k,n) \in H_2(B,\IZ) \oplus (1/2)\IZ \oplus \IZ, 
\]
as stated in Section \ref{intro}.

Let $H$ be a very ample divisor on $B$ and $Z=p^{-1}(H)$. For sufficiently generic $H$ in its linear system, $Z$ is a smooth 
elliptic surface with reduced fibers. Furthermore if $E$ is a vertical pure dimension one sheaf the restriction of $E|_Z$ is a one dimensional sheaf set theoretically suported on a finite union of elliptic fibers. 
Basically $E$ will be said to be generically semistable if 
the restriction $E|_Z=E\otimes_X \CO_Z$ is an 
$\omega|_Z$-semistable pure dimension one sheaf on $Z$ 
for any sufficiently generic very ample divisor $H$ on $B$. 
Technically, this notion requires a more careful definition. 

First note that given any very ample line bundle $L$ on $B$ the projection 
formula yields an  
isomorphism $H^0(X,p^*L)\simeq H^0(B,L)$ since $p_*\CO_X \simeq \CO_B$. Therefore 
the linear system $|L|$ parametrizes simultaneously divisors 
$H\subset B$ as well as vertical divisors $Z=p^{-1}(H)$ in $X$. 
Let $S_{\sf sm}\subset |L|$ denote the open subset parametrizing 
smooth divisors $H_s$ such that $Z_s=p^{-1}(H_s)$ is 
a smooth elliptic surface with reduced fibers.

Since $E$ is vertical of pure dimension two, its scheme theoretic support will be a divisor $D_E$ on $X$ of the form 
\[ 
D_E = \sum_{i=1}^k \ell_ip^{-1}(C_i) 
\]
where $\ell_i \in \IZ$, $\ell_i \geq 1$, and $C_i$ is a reduced  irreducible 
divisor on $B$ for $1\leq i \leq k$. 
Given any very ample line bundle $L$ on $X$ there 
is a nonempty open subset $S_{E,{\sf tr}}\subset S_{\sf sm}$ 
such that the following hold for any closed point $s\in  S_{E, {\sf tr}}$

$(T.1)$ the corresponding divisor $H_s$ intersects each 
$C_i$ transversely at finitely many smooth points of $C_i$, for $1\leq i 
\leq k$ and

$(T.2)$ $H_s$ also intersects the discriminant $\Delta\subset B$ of the map $p:X\to B$ transversely at finitely many smooth points of $\Delta$. This implies that the elliptic fibration
$p|_{Z_s}:Z_s\to H_s$ will be a Weierstrass model with at most nodal fibers. 

\noindent 

Furthermore, according to \cite[Lemma 1.1.13]{huylehn}, 
there exists a second nonempty open subset $S_{E, {\sf pure}}
\subset |L|$ 
such that $E|_{Z_s}$ is a pure dimension one sheaf for any closed
point $s\in S_{E, {\sf pure}}$. 

Before defining generic stability note the following lemma. The proof is straightforward and will be omitted. 

\begin{lemm}\label{vertdimone} 
Let $H$ be a smooth projective curve and $p_Z:Z\to B$ a smooth 
Weierstrass model over $H$. Let $G$ be a pure dimension 
one sheaf on $Z$ with set theoretic support on a finite union 
of elliptic fibers. Let 
$\omega_Z$, $\omega_Z'$ be arbitrary 
 K\"ahler classes on $Z$. Then $G$ is $\omega_Z$-semistable 
if and only if it is $\omega_Z'$ semistable. 
\end{lemm} 

In the situation of Lemma \ref{vertdimone}, the sheaf $G$ will be 
said to be semistable if it is $\omega_Z$-semistable for some 
arbitrary polarization of $Z$. Given a vertical pure dimension two sheaf $E$ and a divisor $Z_s$ corresponding to $s\in S_{E,{\sf tr}} \cap   S_{E, {\sf pure}}$ the sheaf $E|_{Z_s}$ is set theoretically supported on 
a finite union of elliptic fibers. Therefore one can formulate: 

\begin{defi}\label{genstabdef} 
A vertical pure dimension two sheaf $E$ will be called generically 
$\omega$-semistable if and only if for any very ample linear system 
$\Pi=|L|$ on $B$ there exists a nonempty  open 
subset $S_E \subset S_{E,{\sf tr}} \cap   S_{E, {\sf pure}}\subset \Pi$ 
such that the restriction $E|_{Z_s}$ a semistable sheaf on $Z_s$ 
for any closed point $s\in S_E$. 
\end{defi} 

In the remaining part of this section it will be shown that
adiabatic semistability is equivalent to generic semistability 
for vertical semistable pure dimension two sheaves. 
Since the proof is fairly long, it will be divided into 
several shorter steps. 

\begin{lemm}\label{zerotor} 
Let $F$ be an arbitrary pure dimension two sheaf on $X$ and $D\subset X$ a divisor such that $F|_D$ is a one dimensional sheaf on $X$. 
Then ${\mathcal Tor}_k^X(\CO_D,F)=0$ for $k\geq 1$ and 
there is an exact sequence 
\[
0\to F(-D) \to F \to F|_D \to 0. 
\]
where $F(-D) = F\otimes_X \CO_X(-D)$. 
\end{lemm} 

{\it Proof.} This follows immediately from the standard exact sequence 
\[ 
0\to \CO_X(-D) \to \CO_X \to \CO_D \to 0
\]
taking a tensor product by $F$. Under the current assumptions, the sheaf ${\mathcal Tor}_1^X(\CO_D,F)$ is one dimensional, hence 
it must vanish since $F(-D)$ is pure of dimension two. 

\hfill $\Box$

Let $F$ be a nonzero pure dimension two sheaf on $X$ with $\ch_1(F) \in p^*{\rm Pic}(B)$. 
The second Chern class of $F$ is of the form 
\[
\ch_2(F) =\sigma_* (\alpha_F) + k_F f  
\]
where $\alpha_F$ is a curve class on $B$ and $k_F \in (1/2)\IZ$. 
Let 
$H\subset B$ a sufficiently generic very ample divisor on $B$ such that 
$Z=p^{-1}(H)$ is smooth, and 
$F|_Z$ is a pure dimension one sheaf on $Z$. 

\begin{lemm}\label{positivechi} 
Suppose $\chi(F|_Z)>0$. Then 
\[
H \cdot \alpha_F >0.
\] 
\end{lemm}

{\it Proof.} 
Using Lemma \ref{zerotor} and the Rieman-Roch theorem, 
one has
\[
\chi(F|_Z) = \chi(F)-\chi(F(-Z)) = Z\cdot \ch_2(F). 
\]
Then the conclusion follows. 

\hfill $\Box$

\begin{lemm}\label{contralemma} 
Let $\omega = t\Theta - sp^*K_B$, $s,t\in \IR$, $s>t>0$. 
Suppose $E$ is a nonzero adiabatically $\omega$-slope semistable vertical pure dimension 
two sheaf on $X$. 
Let $F\subset E$ be a nonzero proper subsheaf with numerical invariants 
\[
\
\ch_1(F) = p^*C_F, \qquad \ch_2(F) =\sigma_* (\alpha_F) + k_F f  
\]
where  $C_F$ is a nonzero effective divisor class on $B$, 
$\alpha_F$ is an arbitrary divisor class on $B$ and $k_F \in (1/2)\IZ$. 
Then $K_B\cdot \alpha_F \geq 0$. 
\end{lemm}

{\it Proof.} Note that 
\[
\mu_{\omega}(E) = {k \over 
(s-t/2) |K_B\cdot C|}, \qquad 
\mu_{\omega}(F) = 
{(1-s/t) K_B \cdot \alpha_F  + k_F \over 
 (s-t/2)|K_B \cdot C_F|}.
\]
Therefore $F$ destabilizes $E$ for sufficiently small $t>0$ unless 
$ K_B \cdot \alpha_F\geq 0$. 

\hfill $\Box$ 

The proof that adiabatic stability implies generic stability 
uses the same geometric construction as the proof of the Grauert-M\"ulich Theorem 
in \cite[Sect 3.1]{huylehn}. 

Let $L$ be a very ample 
line bundle on $B$, let $V = H^0(B,L)$ and $\Pi=\IP(V)$ denote the 
associated linear system. By convention, $\IP(V) = {\rm Proj}(S^\bullet(V^\vee))$ such that $H^0(\Pi,\CO_\Pi(1)) \simeq V^\vee$. 
Let $\CK$ be the kernel of the evaluation map ${\sf ev}: V\otimes \CO_B \twoheadrightarrow L$, which is a locally free 
sheaf on $B$. 
According to \cite[Sect 3.1]{huylehn}, the total space $\CH$ of the 
projective bundle 
$\IP(\CK)$ parametrizes 
pairs $(H,b)$ with $H\in \Pi$ and $b\in H$ a closed point. 
In more detail, note that the 
evaluation map determines tautologically a section $\theta$ of the line 
bundle $\pi^*L \otimes \pi^*\CO_\Pi(1)$, 
where $\pi:\Pi\times B\to B$ is the canonical projection. Then $\CH$ 
is 
the divisor $\theta =0$ in $\Pi \times B$. In particular there are natural projections
\[
\xymatrix{
\CH \ar[d]_-{\rho} \ar[r]^-{q} & B \\
\Pi }
\]
Moreover, for any closed point $s\in \Pi$
the scheme theoretic inverse image $\rho^{-1}(s)$ is the divisor 
$\theta|_{B_s}=0$ in $B_s = B \times \{s\}\subset B \times \Pi$. 
Given the construction of $\theta$, it follows that the restriction 
$q|_{\rho^{-1}(s)}$ maps $\rho^{-1}(s)$ isomorphically onto 
$H_s$. Let $q_s: \rho^{-1}(s) {\buildrel \sim \over \longto} H_s$ denote the resulting isomorphism. 

For future reference it will be 
useful to provide an explicit construction for the inverse morphism 
$q_s^{-1}: H_s\to \rho^{-1}(s)$. 
By restriction to $H_s$ one obtains an exact sequence 
\[
0 \to  \CK|_{H_s}\to V \otimes \CO_{H_s}  
{\buildrel {\sf ev}_s\over \longto}  L|_{H_s} \to 0
\]
where ${\sf ev}_s= {\sf ev}|_{H_s}$. Let $0\neq z_s\in V$ be a defining section  of $H_s$. Then ${\sf ev}_s(z_s\otimes 1) =0$, hence there is a 
section $y_s \in H^0(H_s,\CK|_{H_s})$ such that the following diagram is 
commutative 
\[
\xymatrix{
0 \ar[r] & \CO_{H_s} \ar[r]^-{1}\ar[d]^-{y_s} & 
\CO_{H_s}\ar[d]^-{f_s} \ar[r] & 0\ar[d] & \\
0 \ar[r] & \CK|_{H_s} \ar[r] &   V \otimes \CO_{H_s}  \ar[r]^-{{\sf ev}_s}
& L|_{H_s} \ar[r] & 0, }
\]
where $f_s(1) = z_s\otimes 1$. Then the snake lemma yields an 
exact sequence 
\[
0\to {\rm Coker}(y_s) \to  {\rm Coker}(f_s)
\to L|_{H_s} \to 0
\]
where ${\rm Coker}(f_s)$ is locally free since $f_s$ is injective on fibers. 
This implies that ${\rm Coker}(y_s)$ is also locally free, hence
$y_s$ is injective on fibers. Therefore $y_s$ determines a section 
$\xi_s:H_s\to \IP(\CK|_{H_s})$. The scheme theoretic 
image of $\xi_s$ coincides tautologically with $\rho^{-1}(s)$ and 
\[
q_s \circ \xi_s = 1_{H_s}. 
\]
Note that $H^0(B, \CK)=0$, hence the section $y_s$ does not 
extend to $B$. However, the following lemma shows that $\xi_s$ 
can be extended to a certain open subset $U\subset B$. 

\begin{lemm}\label{sectlemma} 
There exists an open subscheme $U \subset B$
and a section $\xi: U \to q^{-1}(U)\subset \CH$ 
such that 
\be\label{eq:secteq}
\xi|_{H_s\cap U} = \xi_s|_{H_s \cap U}.
\ee
Furthermore suppose $C\subset B$ is a fixed effective divisor such that 
the set theoretic intersection $C\cap H_s$ 
is a finite set of closed points. Then the open 
subscheme $U$ can be chosen such that $C \subset U$. 
\end{lemm}

{\it Proof.} 
Let $M$ be a very ample line bundle on $B$ such that 
$H^1(B,\CK\otimes_B L^{-1}\otimes_B M)=0$. 
The exact sequence 
\[
0 \to \CK\otimes_B L^{-1}\otimes_B M\to \CK \otimes_B M \to 
(\CK \otimes_B M)|_{H_s} \to 0
\]
yields a surjective map
\be\label{eq:surjmap}
H^0(B, \CK \otimes_B M) \twoheadrightarrow H^0(H_s, (\CK \otimes_B M)|_{H_s}).
\ee
Since $M$ is very ample, $H^0(H_s, M|_{H_s})$ is nontrivial. 
Let $\psi: \CO_{H_s} \to M|_{H_s}$ be a nonzero section of $M|_{H_s}$
and let $U_\psi \subset H_s$ be the complement of the zero divisor 
of $\psi$. Then 
\[
y_s \otimes \psi|_{U_\psi} : \CO_{U_\psi} \to \CK|_{U_\psi} 
\otimes_{U_\psi} M|_{U_\psi}
\]
determines a section of $\IP(\CK|_{U_\psi})$ which coincides with 
$\xi_s|_{U_\psi}$. 
Since the map \eqref{eq:surjmap} is surjective, there exists a 
nonzero
section $y: \CO_B \to \CK\otimes_B M$ such that $y|_{H_s} =y_s\otimes \psi$. Let $I \subset M$ be the image of the morphism $\CK^\vee \to M$ determined by $y$. Then $I \simeq \CI_Y\otimes M$, where 
$\CI_Y$ is the ideal sheaf of a zero dimensional subscheme $Y \subset B$. Let $U\subset B$ be the complement of $Y$. Then $I|_U$ is 
locally free, hence it determines a section $\xi: U \to \IP(\CK)$ 
which agrees with $\xi_s$ over $U_\psi$. 

Suppose $C\subset B$ is a fixed effective divisor which intersects 
$H_s$ at finitely many points. Then for sufficiently ample $M$ the section 
$\psi$ can be chosen such that $U$ contains the set theoretic intersection $C\cap H_s$. Moreover the extension $y$ can be chosen 
such that the support of $Y$ is disjoint from the support of $C$. 

\hfill $\Box$ 

Analogous considerations apply to the linear system $|p^*L|$ on $X$. 
Note that $H^0(X,p^*L) \simeq H^0(B,L)=V$ since 
$p_*p^*L \simeq L$ and there is an exact sequence 
\[
0\to p^*\CK \to V \otimes \CO_X \to p^*L \to 0.
\]
Therefore the total space $\CZ$ of 
$\IP(p^*\CK)$
parametrizes pairs $(Z,x)$ with $Z=p^{-1}(H)$ for some 
$H$ in the linear system $\Pi$, and $x\in Z$ a closed 
point. Note that $\CZ \simeq \CH\times_B X$ and 
there are natural projections 
\[
\xymatrix{
\CZ \ar[d]_-{\rho_Z} \ar[r]^-{q_Z} & X. \\
\Pi }
\]
For any closed point $s\in \Pi$ there is an isomorphism 
 $\rho_Z^{-1}(s)\simeq \rho^{-1}(s)\times_B X$ and 
the restriction 
$q_Z|_{\rho_Z^{-1}(s)}$ maps $\rho_Z^{-1}(s)$ isomorphically 
to $Z_s$. The inverse morphism is given by the section 
$\zeta_s: Z_s \to q_Z^{-1}(Z_s)$, 
\[
\zeta_s = \xi_s \times_{B}\, 1_X.
\]

Now let $\CF= q_Z^*E$ and let $\CF_S=\CF|_{\pi^{-1}(S)}$ for any open subset $S\subset \Pi$. Recall that $S_{\sf sm}\subset |L|$ denotes the open subset parametrizing 
smooth divisors $H_s\subset B$ such that $Z_s=p^{-1}(H_s)$ is 
a smooth elliptic surface with reduced fibers. 
Then one has the following lemma.

\begin{lemm}\label{famHNlemma} 
The following statements hold for any vertical pure dimension 
two sheaf $E$ on $X$. 

$(i)$ There is a nonempty open subset $S_{\sf fl}\subset\Pi$ such that the restrictions $\rho\big|_{\rho^{-1}(S_{\sf fl})}$ and  
$\rho_Z\big|_{\rho_Z^{-1}(S_{\sf fl})}$ are flat and the fibers $\rho^{-1}(s)$, $\rho_Z^{-1}(s)$, 
 $s\in S_{\sf fl}$ are normal irreducible divisors in $X$, $B$ respectively. 

$(ii)$ There is a nonempty open subset $S_{E, {\sf fl}}\subset S_{\sf fl}$ such that $\CF_S$ is flat over $S$ and $E|_{Z_s}$ is pure one 
dimensional for any $s\in S_{E, {\sf fl}}$. 

$(iii)$ There exists a filtration 
\be\label{eq:relativeHN}
0 \subset \CF_1 \subset \cdots \subset \CF_j = \CF.
\ee
of $\CF$ by coherent sheaves on $\CZ$ which restricts to a relative 
Harder-Narasimhan fibration over a suitable nonempty open subset $S_{E, {\sf hn}}\subset S_{E, {\sf fl}}$. 
\end{lemm}

{\it Proof.} The first two statements are completely analogous to 
\cite[Lemma 3.1.1]{huylehn}. For the third statement note that 
\cite[Thm 2.3.2]{huylehn} implies the existence of a filtration of the form 
\eqref{eq:relativeHN} over the open 
subset $\rho_Z^{-1}(S_{E, {\sf fl}})\subset \CZ$. However, this filtration can be 
extended to a filtration of sheaves on $\CZ$ by successive applications 
of \cite[Ex. 5.15.(d)]{har}. 

\hfill $\Box$

Now one can finally prove:

\begin{lemm}\label{genstablemma} 
Let $\omega = t\Theta - sp^*K_B$, $t,s \in \IR$, $0<t<s$, 
and suppose $E$ is an adiabatically $\omega$-slope semistable vertical 
pure dimension two sheaf on $X$. Then $E$ is generically 
semistable.
\end{lemm} 

{\it Proof.} According to Definition \ref{genstabdef} 
one has to prove the existence  of a nonempty open subset 
$S_E \subset S_{E,{\sf tr}} \cap   S_{E, {\sf pure}}$ 
such that the restriction $E|_{Z_s}$ is a semistable pure dimension one sheaf on $Z_s$ 
for any closed point $s\in S_E$. Note that in Lemma \ref{famHNlemma}
one has   $S_{E, {\sf hn}}\subset 
S_{E,{\sf pure}}$ by construction. 
Let $S_{E, {\sf ss}} \subset S_{E, {\sf hn}} \cap S_{E,{\sf tr}}$ be the 
open subset such that $E|_{Z_s}$ is $\omega|_{Z_s}$ semistable 
for all $s\in S_{E, {\sf ss}}$. If $S_{E, {\sf ss}}$ is nonempty, the 
claim follows. 

Suppose $S_{E, {\sf ss}}$ is empty. Let 
$s\in S_{E, {\sf hn}} \cap S_{E,{\sf tr}}$ be a closed point. 
Hence $H_{s}\subset B$ intersects the effective divisor 
$C_E= \sum_{i=1}^k \ell_i C_i$ at finitely many points. 
Let $\xi:U\to \rho^{-1}(U)$ be a section as in 
 Lemma \ref{sectlemma} such that $C_E \subset U$. Then 
relation \eqref{eq:secteq} holds: 
\[ 
\xi|_{H_s\cap U} =\xi_s|_{H_s\cap U}.
\]
Let $Z_{s}= p^{-1}(H_{s})$. Recall that the projection 
$q_Z : \CZ \to X$ maps $\rho_Z^{-1}(s)$ isomorphically to 
$Z_s$ and the inverse morphism 
$\zeta_s: Z_s \to q_Z^{-1}(Z_s)$ is given by
\[ 
\zeta_s = \xi_s \times_{B} 1_X. 
\]
Let $X_U = p^{-1}(U)$. Then 
$\zeta = \xi\times_U 1_{X_U}$ is a section of $q_Z$ over $X_U$ 
such that 
\be\label{eq:secteqB}
\zeta|_{Z_s\cap X_U} = \zeta_s|_{Z_s\cap X_U}.
\ee
Moreover $D_E=p^{-1}(C_E)$ is a subscheme of $X_U$. 

Let $\zeta_E : D_E \to \CZ$ be the restriction of $\zeta$ to $D_E$. 
Let $\varphi:\CF_1 \hookrightarrow \CF$ be the first term in 
the filtration \eqref{eq:relativeHN}. By construction,  
$\varphi|_{\rho_Z^{-1}(s)}$ is injective. Since $D_E$ is 
a subscheme of $X_U$, using equation \eqref{eq:secteqB}, 
one obtains isomorphisms 
\be\label{eq:zetaisomA}
\zeta_E^*\CF_1 \otimes_X \CO_{Z_s} \simeq \zeta^*\CF_1\otimes_{X_U}\CO_{Z_s\cap X_U} 
 \otimes_{X_U} \CO_{D_E} \simeq 
\zeta_s^{*}\big(\CF_1|_{\rho_Z^{-1}(s)}\big)\big|_{Z_s\cap X_U} 
\otimes_{X_U} \CO_{D_E}. 
\ee
However $\CF_1|_{\rho_Z^{-1}(s)}$ is a subsheaf of $\CF|_{\rho_Z^{-1}(s)}\simeq q_{Z,s}^*(E|_{Z_s})$, where $q_{Z,s} : \rho_Z^{-1}(s) \to Z_s$ 
denotes the natural projection. Since $\zeta_s: Z_s \to 
\rho_Z^{-1}(s)$ is an isomorphism and 
$q_{Z,s}\circ \zeta_s = 1_{Z_s}$, 
it follows that $\zeta_s^{*}\big(\CF_1|_{\rho_Z^{-1}(s)}\big)$ is a 
subsheaf of $E|_{Z_s}$. In particular it is scheme theoretically supported on $D_E$, and equation \eqref{eq:zetaisomA} yields 
an isomorphism 
\be\label{eq:zetaisomB}
\zeta_E^*\CF_1 \otimes_X \CO_{Z_s}\simeq \zeta_s^{*}\big(\CF_1|_{\rho_Z^{-1}(s)}\big).
\ee
Now let 
$f= \zeta_E^*\varphi: \zeta_E^* \CF_1 \to \zeta_E^*\CF\simeq E$ and 
let $F={\rm Im}(f) \subset E$. Since $\zeta_s^{*}\big(\CF_1|_{\rho_Z^{-1}(s)}\big)$ is a 
subsheaf of $E|_{Z_s}$, equation \eqref{eq:zetaisomB} 
implies that 
\[
F|_{Z_s} \simeq \zeta_E^*(\CF_1|_{\rho_Z^{-1}(s)})
\]
is also a subsheaf of $E|_{Z_s}$. By construction this is the first term in the Harder-Narasimhan filtration of $E|_{Z_s}$. 
Since $E$ is vertical, Lemma \ref{zerotor} implies that 
$\chi(E|_{Z_{s}})=0$. Thefore, as $E|_{Z_s}$ is not semistable by assumption, one must have $\chi(F|_{Z_{s}})> 0$. 
Then Lemma \ref{positivechi} implies that $c_1(L) \cdot \alpha_F >0$, 
where $\alpha_F$ is the horizontal part of $\ch_2(F)$ as in  loc.cit. 

Let $\omega_B$ be the dualizing sheaf of $B$. Applying the above construction to $L = \omega_B^{-m}$, for sufficiently large  $m\geq 1$,  one is then led to a contradiction with Lemma \ref{contralemma} since $E$ is assumed 
adiabatically $\omega$-slope semistable. 

\hfill $\Box$

Lemma \ref{genstablemma} admits the following converse. 

\begin{lemm}\label{adiabstablemma} 
Suppose $\omega = t\Theta -s K_B$ with $s,t\in \IR$, 
$s >t>0$. 
Let $E$ be an $\omega$-semistable vertical pure dimension 
two sheaf on $X$ with numerical invariants 
\[ 
\ch_1(E) = p^*C, \qquad \ch_2(E) = kf
\]
where $C$ is a nonzero
 effective divisor class on $B$ and $k\in (1/2)\IZ$. 
Suppose $E$ is generically semistable. Then $E$ is 
adiabatically $\omega$-semistable. 
\end{lemm} 

{\it Proof.} Let $0\neq F \subset E$ be a proper pure dimension two 
subsheaf of $E$ such that $G=E/F$ is also pure dimension two. 
Let $H$ be a sufficiently generic very ample divisor on $B$ as in Definition
\ref{genstabdef} and $Z=p^{-1}(H)$. 
Lemma \ref{zerotor} shows that ${\mathcal Tor}_1^X(G,\CO_Z)=0$, 
hence there is an exact sequence 
\[
0\to F|_Z \to E|_Z \to G|_Z \to 0. 
\]
Let $\ch_1(F)=p^*C_F$ and $\ch_2(F) =\sigma_*(\alpha_F) + k_F f$ 
with $C_F,\alpha_F$ divisor classes on $B$, $C_F$ nonzero, effective, 
and $k_F \in (1/2)\IZ$. 
Since $E|_Z$ is semistable by assumption, and $\chi(E|_Z)=0$, it follows 
that 
\[
\chi(F|_Z) = Z\cdot \ch_2(F) = H \cdot \alpha_F \leq 0.
\]
In particular, for $H$ in the linear system $|-K_B|$, 
\be\label{eq:negproduct}
K_B \cdot \alpha_F \geq 0.
\ee

Let $\omega'= t'\Theta -s p^*K_B$ with $0<t'\leq t$. Then 
\[
\mu_{\omega'}(F) = -{s\over t'(s-t'/2)} {K_B\cdot \alpha_F\over 
|K_B \cdot C_F|} + {1\over s-t'/2} {k_F+ K_B \cdot \alpha_F\over |K_B \cdot C_F|}
\]
and 
\[
\mu_{\omega'}(E) ={1\over s-t'/2}  {k\over 
|K_B\cdot C|}.
\]
Since $\mu_{\omega}(F) \leq \mu_{\omega}(E)$, 
one finds 
\[
{k_F+K_B \cdot \alpha_F \over |K_B\cdot C_F|} \leq {k\over |K_B\cdot C|} + {s\over t} {K_B\cdot \alpha_F \over |K_B\cdot C_F|}.
\]
Using \eqref{eq:negproduct}, this  implies that 
\[
\mu_{\omega'}(F) -  \mu_{\omega'}(E) \leq 
{s(t'-t)\over tt'(s-t'/2)} {K_B\cdot \alpha_F \over |K_B \cdot C_F|}\leq 0
\]
for any $0<t'<t$. 
Moreover equality holds for some $0<t'<t$ 
if and only if $K_B \cdot \alpha_F=0$. 
If this is the case, $\omega$-stability implies that 
\[
\nu_\omega (F) \leq \nu_\omega(E) 
\]
which is equivalent to 
\[ 
{\chi(F)\over |K_B\cdot C_F|} \leq {\chi(E) \over |K_B\cdot C|}.
\]
This implies that $\nu_{\omega'}(F) \leq \nu_{\omega'}(E)$. 
Therefore $E$ is $\omega'$-semistable. 

\hfill $\Box$

Let $\CM_{\omega}(X,\gamma)$ denote the 
moduli stack of $\omega$-semistable pure dimension 
two sheaves $E$ with numerical invariants 
 $\gamma=(C,k,n)\in H_2(B,\IZ)\oplus (1/2)\IZ\oplus \IZ$.
Let $\CM^{\sf ad}_{\omega}(X,\gamma)$ denote the substack 
of adiabatically semistable objects. 
To conclude this section it will be shown that 
$\CM^{\sf ad}_{\omega}(X,\gamma)$ is 
an open substack of $\CM_{\omega}(X,\gamma)$
for any discrete invariants $\gamma$ and  for any K\"ahler class 
$\omega = t \Theta - sp^*K_B$, $s,t, \in \IR$, $s>t>0$. 
For any $0< t'<t<s$ let $\omega_{t'} = t'\Theta -s p^*K_B$.
Then one has:

\begin{lemm}\label{ttprimelemma} 
Suppose $E$ is a vertical $(\omega_t,\beta)$-semistable 
sheaf with discrete invariants $\gamma=(C,k,n)$, 
$C\neq 0$, which is  not 
$\omega_{t'}$-semistable for some $0<t'<t$.  Then 
$E$ is not $\omega_{t''}$-semistable for any $0< t'' <t'$. 
\end{lemm}

{\it Proof}. Let $F\subset E$ be a destabilizing proper non-zero subsheaf 
with respect to $\omega_{t'}$-stability. This means that 
\be\label{eq:tprimedestab}
\mu_{\omega_{t'}}(F) \geq  \mu_{\omega_{t'}}(E), 
\ee
and, if equality holds, $\nu_{\omega_{t'}}(F) > \nu_{\omega_{t'}}(E)$. 
At the same time, note that 
$\mu_{\omega_{t}}(F) \leq  \mu_{\omega_{t}}(E)$.
As in the proof of Lemma \ref{adiabstablemma}, let 
\[
\ch_1(F) = p^*(C_F), \qquad \ch_2(E') = \sigma_*(\alpha_F) + k_Ff
\]
where $C_F$ is a nonzero effective curve class on $B$. 
Then the same computation as in loc.cit. shows that  
\[
\mu_{\omega_{t'}}(F) -  \mu_{\omega_{t'}}(E) \leq 
{s(t'-t)\over tt'(s-t'/2)} {K_B\cdot \alpha_F\over |K_B\cdot C_F|}.
\]
Therefore inequality \eqref{eq:tprimedestab} implies that $K_B \cdot \alpha_F 
\leq 0$. Moreover, if $K_B \cdot \alpha_F=0$, equality must hold in 
\eqref{eq:tprimedestab}.

Now suppose $E$ is $\omega_{t''}$-semistable for some 
$0< t'' < t'$. Then 
\be\label{eq:tprimeprime}
\mu_{\omega_{t''}}(F) \leq  \mu_{\omega_{t''}}(E).
\ee
and, if equality holds, $\nu_{\omega_{t''}}(F) < \nu_{\omega_{t''}}(E)$. 
However inequality \eqref{eq:tprimedestab} yields 
\[ 
{k_F+K_B \cdot \alpha_F\over |K_B\cdot C_F|} \geq  {k\over |K_B\cdot C|} + {s\over t'} {K_B\cdot \alpha_F \over |K_B \cdot C_F|}.
\]
Therefore 
\[
\mu_{\omega_{t''}}(F)-\mu_{\omega_{t''}}(E) \geq 
{s(t''-t')\over t't''(s-t''/2)} {K_B\cdot \alpha_F \over |K_B \cdot C_F|} \geq 0. 
\]
This implies that $K_B\cdot \alpha_F=0$, hence equality must hold in 
\eqref{eq:tprimedestab} and \eqref{eq:tprimeprime}. However, in this case, $\nu_{\omega_{t'}(F)} > 
\nu_{\omega_{t'}}(E)$, which is equivalent to 
\[
{\chi(F)\over |K_B\cdot C_F|} > {\chi(E)\over |K_B\cdot C_F|}.
\]
This further implies $\nu_{\omega_{t''}(F)} > 
\nu_{\omega_{t''}}(E)$, leading to a contradiction. 

\hfill $\Box$ 

In order to formulate the last result of this section, let 
$M_\omega(X,\gamma)$ denote the coarse moduli scheme 
parameterizing $S$-equivalence classes of $\omega$-semistable 
sheaves on $X$. Note that according to \cite[Ex.8.7]{goodmoduli}, 
$M_\omega(X,\gamma)$ is a good coarse moduli space for the moduli 
stack $\CM_{\omega}(X,\gamma)$. This means that there is a 
 morphism $\varrho: \CM_{\omega}(X,\gamma)
\to M_{\omega}(X,\gamma)$ satisfying the properties listed in \cite[Thm. 4.16]{goodmoduli}. Let $M_\omega^{\sf ad}(X,\gamma)\subset M_\omega(X,\gamma)$ be 
the scheme theoretic image $\varrho(\CM_\omega^{\sf ad}(X,\gamma))$.

\begin{lemm}\label{adiabopen}
For any K\"ahler class $\omega = t\Theta-sp^*K_B$ with 
 $s>t>0$, and 
for any discrete invariants $\gamma$, the subscheme $M_\omega^{\sf ad}(X,\gamma)$ is open in $M_\omega(X,\gamma)$, and the substack 
$\CM^{\sf ad}_{\omega}(X,\gamma)$ is open in $\CM_{\omega}(X,\gamma)$.
\end{lemm}

{\it Proof.} For any $0<t'<t$, let $\CN_{t'}(\gamma)$ be the substack of $\omega$-semistable 
vertical sheaves which are not $\omega_{t'}$-semistable.  
Note that this is an closed substack of $\CM_{\omega}(X,\gamma)$
since $\omega_{t'}$-semistability is an open condition in flat families. 
Lemma \ref{ttprimelemma} shows that  $\CN_{t'}(\gamma)\subset 
\CN_{t''}(\gamma)$ for any $0<t''<t'<t$. According to 
\cite[Thm 4.16.(i)]{goodmoduli}, 
the morphism $\varrho$ is universally closed. Therefore 
the scheme theoretic image $\varrho(\CN_{t'}(\gamma))$ is a 
closed subscheme $N_{t'}(\gamma)\subset  M_\omega(X,\gamma)$. 
Moreover and \cite[Thm.4.16.(iii)]{goodmoduli} 
implies that $N_{t'}(\gamma)\subseteq
N_{t''}(\gamma)$ for any $0<t''<t'<t$. Since $M_\omega(X,\gamma)$ is noetherian, 
it follows that the union $N_\omega(\gamma)=\cup_{0<t'<t} N_{t'}(\gamma)$ 
must be a closed subscheme of $M_\omega(X,\gamma)$. Therefore its inverse 
image $\CN_\omega(\gamma)=\varrho^{-1}(N_\omega(\gamma))$ is a closed substack of $\CM_\omega(X,\gamma)$. 
To conclude the proof note that 
$\CM_\omega^{\sf ad}(X,\gamma)$ is the complement of 
$\CN_\omega(\gamma)$ according to Lemma \ref{ttprimelemma}.

\hfill $\Box$

\section{Fourier-Mukai transform and stability}\label{FMvertical}

This section contains the detailed proof of Theorem \ref{main}. 
Since the proof is fairly long and complicated, it will be divided 
into subsections. The first subsection reviews the basic proeprties of 
the relative Fourier-Mukai transform on elliptic fibrations. 

\subsection{Basics of  Fourier-Mukai transform}\label{FMbasics}

The main references for this section will be \cite{FM_thesis, 
FM_elliptic, FM_K3_elliptic} and the review article \cite{FM_string}. 
Let $X$ be a smooth generic elliptic Weierstrass model over a smooth Fano variety $B$. In particular all singular elliptic fibers are either 
nodal of cuspidal. In this subsection $X$ will be assumed of dimension $n\in \{2,3\}$ and not necessarily Calabi-Yau. 
Let ${\hat X}$ be the Altman-Kleiman compactification of the degree zero relative Jacobian of $X$
and ${\hat p}:{\hat X}\to B$ the natural projection.  
This is a fine relative moduli space for rank one degree zero 
torsion free sheaves on the fibers of $p$, hence there is  
a (non-unique)  
universal 
rank one torsion free sheaf $\calP$ on $\hX\times_B {X}$. 
There is also a canonical morphism $\theta : X\to {\hat X}$ 
mapping  a closed point $x\in X$ to 
$\CI_x \otimes \CO_{X_{p(x)}}(\sigma(p(x)))$, where $\CI_x\subset 
\CO_{X_{p(x)}}$ 
is the ideal sheaf of $\{x\}\subset X_{p(x)}$, and 
$\sigma: B \to X$ is the canonical section. 
Under the current assumptions $\theta$ is an isomorphism. 
Hence ${\hat p}:{\hat X}\to B$ 
is a smooth Weierstrass model with a canonical section ${\hat \sigma}: 
B \to {\hat X}$. 

Note that $\calP$ is flat 
over ${\hat X}$ and also flat over $X$ according to 
\cite[Lemma 8.4]{FM_K3_elliptic}. 
Extending $\calP$ by zero to $\hX\times X$, let 
\be\label{eq:dualkernel}
\CQ = {\rm R}{\mathcal Hom}_{\hX\times X} (\calP, 
\pi_X^*\omega_X)[n-1] 
\ee
where $\pi_X : \hX \times X \to X$ is the canonical projection
and $\omega_X$ is the dualizing sheaf of $X$. 
Then \cite[Lemma 8.4]{FM_K3_elliptic} proves that 
$\CQ$ is a sheaf on $\hX \times X$ which is flat over both 
$\hX$ and $X$. Moreover, $\CQ$ is pure and scheme theoretically 
supported on $\hX\times_B X$. 

Now consider the commutative diagram 
\be\label{FMdiagA} 
\xymatrix{ 
{ \hX} \times_B X \ar[r]^-{ \rho} \ar[d]^-{\hat \rho} \ar[dr]^-{q}  & {X} \ar[d]^-{p}
\\ 
\hX \ar[r]^-{\hat p} & B.\\ }
\ee
and define the Fourier-Mukai functors $\Phi:D^b(\hX)\to D^b(X)$, 
\be\label{eq:FMa} 
\Phi(\hE) = R{\rho}_*(L{\hat \rho}^*(\hE) {\buildrel L\over \otimes}\calP) 
\ee
and ${\widehat \Phi}: D^b(X) \to D^b(\hX)$, 
\be\label{eq:FMb} 
{\widehat \Phi}(E) = R{\hat \rho}_*(L{\rho}^*(E)
{\buildrel L\over \otimes} {\mathcal Q}).
\ee
Theorem \cite[Thm 1.2]{FM_K3_elliptic} proves the
following relations:  
\be\label{eq:InverseFM} 
{\widehat \Phi}\circ \Phi \simeq {\rm Id}_{D^b(\hX)}[-1]\qquad 
{\Phi}\circ {\widehat \Phi} \simeq {\rm Id}_{D^b(X)}[-1].
\ee


For any object $E$ in $D^b(X)$ let ${\widehat \Phi}^{i}(E)$ denote the $i$-th cohomology sheaf of ${\widehat \Phi}(E)$.
Since ${\CQ}$ is flat over $X$, the base change theorem 
implies that $\hPhi^i(E)$ is nonzero only for $i\in \{0,1\}$.  
A sheaf $E$ on $X$ is called ${\widehat \Phi}-WIT_i$ if 
${\widehat \Phi}^j(\hE)=0$ for all $j\neq i$. The same 
applies to sheaves on $\hX$ with respect to the inverse  
functor $\Phi$.  

For any closed point ${\hat x}\in \hX$ let 
$\iota_{\hat x}: {\hat x}\times X 
\hookrightarrow \hX\times X$ denote the canonical embedding
and $\calP_{\hat x} = \iota_{\hat x}^*\calP$, $\CQ= \iota_{\hat x}^*\CQ$. Note that $\calP_{\hat x}$ is isomorphic to the extension 
by zero of a rank one torsion free sheaf on the elliptic fiber 
$X_{{\hat p}({\hat x})}$. Since $\calP, \CQ$ are flat over $\hX$, 
\cite[Lemma 3.1.1]{FM_thesis} implies that 
\[
{\rm L}_k \iota_{\hat x}^* {\cal P} =0, \qquad 
{\rm L}_k\iota_{\hat x}^* \CQ=0
\]
for all $k>0$. Then, using \cite[Prop. III.8.8]{resdual}, 
relation \eqref{eq:dualkernel} yields the isomorphism
\be\label{eq:dualityrelA}
\CQ_{\hat x} \simeq  {\rm R}{\mathcal Hom}_{X}({\cal P}_{\hat x}, \omega_X)[n-1]
\ee
in $D^b(X)$. This implies that $\CQ_{\hat x}$ is a pure dimension 
one sheaf on $X$ with scheme theoretic support on $X_{{\hat p}({\hat x})}$. Taking a further derived dual, one also has 
\be\label{eq:dualityrelAB}
\calP_{\hat x} \simeq  {\rm R}{\mathcal Hom}_{X}(\CQ_{\hat x}, \omega_X)[n-1]
\ee
Analogous results hold 
for the fibers of ${\rho}$. 

Next note the following lemma, which is a simple consequence of the definitions. 
\begin{lemm}\label{FMpoint} 
$(i)$  For any closed point ${\hat x}\in \hX$ the skyscraper sheaf $\CO_{\hat x}$ is $\Phi-WIT_0$ and 
\[ 
\Phi^0(\CO_{\hat x}) \simeq {\mathcal P}_{\hat x}. 
\]

$(ii)$ For any closed point ${ x}\in X$
 the $\CO_{\hX}$-module ${\mathcal P}_{x}$ is $\Phi-WIT_1$ and 
\[
\Phi^1({\mathcal P}_{x}) \simeq \CO_{x}. 
\]

$(iii)$ Analogous results hold for closed points $x\in X$ relative to 
$\hPhi$. 
\end{lemm}

Further results needed in the following include 
\cite[Lemma 9.2]{FM_K3_elliptic} and \cite[Lemma 9.3]{FM_K3_elliptic} 
which will be reproduced below for convenience. 

\begin{lemm}\label{onezero} 
Let $\hE$ be a sheaf on $\hX$. Then $\Phi^i(\hE)$ is $\hPhi-WIT_{1-i}$
for $i\in \{0,1\}$ and there is a short exact sequence 
\[ 
0 \to \hPhi^1(\Phi^0(\hE)) \to \hE \to \hPhi^0(\Phi^1(\hE))\to 0. 
\]
An analogous statement holds of sheaves $E$ on $X$ with 
$\Phi$ and $\hPhi$ reversed. 
\end{lemm} 

\begin{lemm}\label{WITO} 
A sheaf $\hF$ on $\hX$ is $\Phi-WIT_0$ if and only if 
${\rm Hom}_\hX(\hF, \CQ_x) =0$ for all $x\in X$.
\end{lemm}

Now suppose $X$ 
is a Calabi-Yau threefold. 
Choosing the normalization of \cite{FM_string} let $\calP$ be 
given by 
\[
\calP= \CI_{\Delta} \otimes \rho^*\CO_X(\Theta) \otimes 
{\hat \rho}^* \CO_\hX({\widehat \Theta}) \otimes q^{*} \omega_B^{-1} 
\]
where $\CI_\Delta$ is the ideal sheaf of the diagonal 
$\Delta \subset X\times_B \hX\simeq X \times_B X$, $\omega_B$ is 
the dualizing sheaf of $B$, and $\Theta\subset X$, ${\widehat \Theta}\subset \hX$ are the canonical sections. This particular choice 
for $\calP$ will be used throughout the remaining part of the paper. 
Then note that equations (17) and (18) in \cite[Sect 5.3]{FM_string} 
yield the following formulas for the Chern characters
of the Fourier-Mukai transform of vertical sheaves.

Let $\hF$ be a pure dimension one sheaf on $\hX$ with
\[ 
\ch_2(\hF) = {\hat \sigma}_*(C) + m {\hat f}, \qquad 
\ch_3(\hF) = l \ch_3(\CO_{\hat x}), \qquad m,l\in \IZ. 
\]
Then 
\be\label{eq:FMchernA} 
\bal 
\ch_0(\Phi(\hF)) = 0, & \qquad \ch_1(\Phi(\hF)) = p^*C, \qquad 
\ch_2(\Phi(\hF)) = \big(l + K_B \cdot C/2\big) f \\
& \ch_3(\Phi(\hF)) = -m \ch_3(\CO_{x}).\\
\eal 
\ee
with $x\in X$ an arbitrary closed point. 
Conversely, let $E$ be a 
vertical pure dimension two sheaf on $X$ with 
\[
\ch_1(E) = p^*C, \qquad \ch_2(E) = k f, 
\qquad \ch_3(E) = -n \ch_3(\CO_x) 
\]
where $C$ is an effective  curve class on $B$ 
and $k\in (1/2)\IZ$, $n \in \IZ$, $k\equiv K_B\cdot C/2$ mod $\IZ$. 
Then 
\be\label{eq:FMchernB} 
\bal 
\ch_0(\hPhi(E)) =0, & \qquad \ch_1(\hPhi(E)) =0, \qquad  \ch_2(\hPhi(E)) = - {\hat \sigma}_*C - n {\hat f}  \\
& \ch_3(\hPhi(E)) = \big(-k +  K_B\cdot C/2\big) \ch_3(\CO_{\hat x}).\\
\eal 
\ee

\subsection{From sheaves on $X$ to sheaves on $\hX$} 

\begin{lemm}\label{Reslemma} 
Let $E$ be a vertical pure dimension two 
 sheaf on $X$. 
Let $\hU\subset \hX$ be an arbitrary affine open subset. 
Then $\hPhi(E)|_\hU$ is quasi-isomorphic to a three term complex 
of coherent locally free $\CO_\hU$-modules 
\[
0\to W_{-1} {\buildrel \phi_0\over \longto} 
 W_0 {\buildrel \phi_1 \over \longto}  W_1 \to 0
\]
where the degree of $W_{i}$ is $i$ for $-1\leq i \leq 1$ 
and $\phi_0$ is injective. 
\end{lemm} 

{\it Proof.} Since $E$ is pure dimension 
two, it has a locally free resolution 
\[
V_{-2}\to V_{-1}\to V_0 
\]
on $X$, where $V_{-i}$ is in degree $-i$ for $0\leq i\leq 2$. 
Since $\rho$ is flat and $\CQ$ is flat over $X$, 
$L{\rho}^*(E)
{\buildrel L\over \otimes} {\mathcal Q}$ is isomorphic to the 
complex 
\[
 \rho^*V_{-2} \otimes \CQ \to \rho^*V_{-1} \otimes \CQ
\to \rho^*V_0 \otimes \CQ
\]
in $D^b(\hX\times_B X)$. Let $\CV$ denote the above complex and note that each
term of this complex is flat over $\hX$. 

Given any affine open subscheme $\hU\subset \hX$ let 
${\hat \rho}_\hU$ denote the restriction 
of ${\hat \rho}$ to ${\hat \rho}^{-1}(\hU)$. 
Then $\hPhi(E)|_\hU$ is given by 
${\rm R}{\hat \rho}_{\hU*}(\CV|_{{\hat \rho}^{-1}(\hU)})$. 
According to \cite[Thm. 6.10.5]{EGAIIIb} and 
\cite[Remark 6.10.6]{EGAIIIb}, or 
\cite[Sect 5, page 46]{abelian_var},  ${\rm R}{\hat \rho}_{\hU*}(\CV|_{{\hat \rho}^{-1}(\hU)})$ is quasi-isomorphic to a finite complex  
$W_{\bullet}$ of locally free $\CO_\hU$-modules. 
Moreover, for any point ${\hat x}\in \hU$, the cohomology group 
$H^i(W_\bullet|_{\hat x})$ is isomorphic to the hypercohomology group 
$\IH^i(\CV|_{{\hat \rho}^{-1}({\hat x})})$ for all values of $i$. 

Next note  that $W_\bullet$ can be 
truncated to a three term locally free complex of amplitude 
$[-1,\ 1]$. By general properties of the Fourier-Mukai transform, $\hPhi(E)|_\hU$ has nontrivial 
cohomology sheaves only in degrees $0,1$ hence one can truncate $W$ to a locally free complex 
\[
\cdots \to W_{-1} \to W_{0} \to W_1 \to 0
\]
where $W_i$ is in degree $i$ for all $i\leq 1$. Recall that the cokernel of an injective morphism $f_i: W_i \to W_{i-1}$ of locally free sheaves is locally free 
if and only if $f_i$ is injective on fibers. 
Then the claim will follow 
if one shows that 
\[ 
\IH^{-i}(\CV|_{{\hat \rho}^{-1}({\hat x})})=0
\]
for all $i\geq 2$. In order to prove this, 
note that 
the cohomology sheaf 
$\CH^{-i}(\CV|_{{\hat \rho}^{-1}({\hat x})})$ is isomorphic to the 
local Tor sheaf ${{\CT}or}_i^X(E, \CQ_{\hat x})$.
 Then relation \eqref{eq:dualityrelA} yields isomorphisms 
\be\label{eq:sheafcohA}
\CH^{-i}(C|_{{\hat \rho}^{-1}({\hat x})})\simeq
{\mathcal Ext}^{2-i}_X \big({\calP}_{\hat x}, E) 
\ee
for all $ i\in \IZ$.
In particular $\CH^{-i}(C|_{{\hat \rho}^{-1}({\hat x})})=0$, 
$i\geq 3$ for degree reasons, and 
$\CH^{-2}(C|_{{\hat \rho}^{-1}({\hat x})})=0$ since $E$ is pure dimension two. Then the required 
vanishing result follows from the hypercohomology spectral 
sequence since the remaining cohomology sheaves 
of $\CV|_{{\hat \rho}^{-1}({\hat x})}$ are set theoretically supported 
in dimension one. 
In conclusion $\hPhi(E)|_U$ is quasi-isomorphic 
to a complex of the form 
\[
0\to W_{-1} {\buildrel \phi_0\over \longto} 
 W_0 {\buildrel \phi_1 \over \longto}  W_1 \to 0
\]
where $\phi_0$ is injective. 

\hfill $\Box$.

\begin{lemm}\label{genlemmaA} 
Let $H$ be a smooth projective curve and $Z$ be a smooth 
Weiestrass model over $H$ with at most nodal fibers. Let $F$ be a coherent sheaf on $Z$ 
with set theoretic support on a reduced fiber $Z_b$, for $b\in H$  
a closed point. Suppose $\chi(F)=0$ and $F$ is stable with respect to an arbitrary polarization $\omega_Z$. Then $F$ is the extension 
by zero of a rank one torsion free sheaf $G$ on $Z_b$ with $\chi(G)=0$. 
\end{lemm} 

{\it Proof.} 
According to \cite[Thm 1.1]{Vb_elliptic}, any stable torsion free sheaf $G$ on $Z_b$ with $\chi(G)=0$ must
have rank one. Therefore it suffices to prove that $F$ is scheme theoretically supported on $Z_b$. 

Since $F$ is stable, it must be 
pure of dimension one. Hence it is scheme theoretically supported
on a nonreduced divisor $kZ_b$ on $Z$ for some $k\in \IZ$, 
$k>0$. Consider the morphism $F {\buildrel \zeta\over \longto} 
F\otimes_Z \CO_Z(Z_b)$, where $\zeta \in \CO_Z(Z_b)$ is a defining section.  Note that $\CO_{Z}(Z_b)\simeq \CO_Z(Z_{b'})$ for 
any point $b'\in H\setminus \{b\}$. Pick any such point $b'$ and 
let $\zeta'\in\CO_Z(Z_b)$ be its defining section. Obviously 
$\zeta'$ is nonzero on $Z \setminus Z_{b'}$, hence its yields 
an isomorphism $F\otimes_H \CO_Z(Z_b)\simeq F$. 
Since $F$ is assumed stable it follows that $F {\buildrel \zeta\over \longto} 
F\otimes_H \CO_Z(Z_b)$ must be either identically zero or an isomorphism. However note that $F {\buildrel \zeta^k\over \longto} 
F\otimes_Z \CO_Z(kZ_b)$ must be identically zero since $F$ is 
scheme theoretically supported on $kZ_b$. Therefore  $F {\buildrel \zeta\over \longto} F\otimes_Z \CO_Z(Z_b)$ cannot be an isomorphism, which implies 
that it must be identically zero. In conclusion $F$ is scheme 
theoretically supported on the reduced fiber $Z_b$, hence it must be isomorphic to the extension by zero of a stable 
sheaf $G$ on $Z_b$. 

\hfill $\Box$ 

Let $E$ be a vertical pure dimension two 
 sheaf on $X$ scheme theoretically supported on a divisor 
\be\label{eq:schsupp}
 D_E=\sum_{i=1}^k \ell_i p^{-1}(C_i)
\ee
for some reduced irreducible effective divisors $C_i$ in $B$. 
Using the notation introduced above Definition \ref{genstabdef}, 
let $H$ be a very ample divisor in $B$ corresponding to a closed 
point $s\in S_{E, {\sf tr}}\cap S_{E, {\sf pure}} \subset |H|$. 
Therefore $Z=p^{-1}(H)$ is a smooth elliptic surface with finitely many nodal fibers which intersects each component 
$p^{-1}(C_i)$ transversely 
along a finite collection of elliptic fibers. 

Next note that ${\widehat Z} = {\hat p}^{-1}(H)\subset {\hX}$ is a smooth 
elliptic surface isomorphic to $Z$ over $H$.  Moreover
$\hZ\times_H Z = (\hX \times_B H) \times_H Z = \hX\times_H Z$ is the inverse image $\rho^{-1}(Z)$ under the projection 
$\rho: \hX \times_B X \to X$. In particular $\hZ\times_H Z$ is a 
closed subscheme of $\hX\times_B X$ and $\rho^*\CO_Z \simeq 
\CO_{\hZ\times_H Z}$. 
Let 
\be\label{eq:canmorphism}
\jmath: {\widehat Z}\times_H Z \to \hX\times_B \hX
\ee
denote the canonical closed embedding and let 
$\hPhi_Z: D^b(Z) \to D^b(\hZ)$ be the Fourier-Mukai functor 
with kernel ${\rm L}\jmath^*\CQ$.

\begin{lemm}\label{genvertical} 
Suppose $E$ is a nonzero vertical pure dimension two sheaf 
with scheme theoretic support \eqref{eq:schsupp} and let 
$Z=p^{-1}(H)\subset X$ be a vertical divisor as above. 
Then there is an isomorphism
\be\label{degoneisomA} 
\hPhi^1(E) \otimes_X \CO_\hZ \simeq \hPhi_Z^1(E\otimes_X \CO_Z) 
\ee 
\end{lemm} 

{\it Proof}. 
Since $\rho: \hX\times_B X \to X$ is flat, there is an exact sequence 
\[ 
0\to \rho^* \CO_X(-Z) \to \rho^*\CO_X \to \rho^*\CO_Z\to 0
\]
where $\rho^*\CO_X \simeq \CO_{\hX\times_B X}$ and 
$\rho^*\CO_Z\simeq \CO_{\hZ\times_HZ}$. Hence this is a two 
term locally free resolution of $\CO_{\hZ\times_HZ}$. 
Since $\CQ$ is flat over $X$, this sequence will remain exact when 
one takes a tensor product with $\CQ$. Therefore $L\jmath^*\CQ$ is quasi-isomorphic to 
$\CQ|_{{\widehat Z}\times_H Z}$.

Since the Fourier-Mukai transform is compatible with base change there is an isomorphism 
\[ 
\hPhi_Z({\rm L}\iota_Z^*(E) ) \simeq L\iota_{{\hZ}}^*\hPhi(E)
\]
in $D^b(\hZ)$, 
where $\iota_{\hZ}: {\hZ}\to  \hX$ is the natural closed 
embedding. However, Lemma \ref{zerotor} yields an 
isomorphism ${\rm L}\iota_Z^*E \simeq E\otimes_X \CO_Z$ 
in $D^b(Z)$, hence one obtains 
\[ 
\hPhi_Z(E\otimes_X \CO_Z) \simeq L\iota_{{\widehat Z}}^*\hPhi(E).
\]
Since $\hPhi(E)$ has cohomology only in degrees 0 and 1, the base change theorem \cite[Thm. 7.7.5]{EGAIIIb} implies that 
\be\label{eq:HrestrictionB}
\hPhi^1_Z(E\otimes_X \CO_Z) \simeq \iota_{\widehat Z}^*\hPhi^1(E) \simeq 
\hPhi^1(E)\otimes_{\hX} \CO_{\widehat Z}. 
\ee

\hfill $\Box$ 

\begin{lemm}\label{FMpuredimtwo}
Let $E$ be a nonzero generically semistable vertical pure dimension two sheaf on $X$ as 
in Definition 
\ref{genstabdef}. 
Then $\hPhi^0(E)=0$ and $\hPhi^1(E)$  is a pure dimension one 
sheaf on $\hX$. 
\end{lemm} 

{\it Proof.} Using the notation of Lemma \ref{Reslemma} 
it suffices the prove that the complex $W_{\bullet}$ is exact in degree 0
for any open affine subset $\hU\subset X$. 
 Under the current assumptions the scheme theoretic 
support of $E$ is of the form 
\eqref{eq:schsupp}. Note that the first Chern character 
of $E$ is of the form 
\be\label{eq:firstchernA}
\ch_1(E) =\sum_{i=1}^k r_i p^*(C_i)
\ee
for some integers $r_i\in \IZ$, $r_i\geq 1$, $1\leq i\leq k$. 

Let $H$ be a smooth very ample divisor on $B$ 
satisfying the genericity conditions in Definition \ref{genstabdef}. 
 In particular $H$ intersects 
each $C_i$ transversely at $n_i\geq 1$ finitely many smooth points $b_{i,j}$ on $C_i$, where $1\leq j \leq n_i$. The inverse image 
$Z=p^{-1}(H)$ is a smooth Weierstrass model over $H$ with at most 
nodal fibers. 
Let $F=E\otimes_X \CO_Z$. 
By assumption, 
$F$ is an $\omega|_{Z}$- semistable sheaf on $Z$ set theoretically 
supported on a finite union of elliptic fibers. 
According to  Lemma \ref{zerotor}, 
there is an exact sequence 
\[
0\to 
E(-Z) \to E \to F\to 0
\]
of sheaves on $X$ which yields $\chi(F) =0$ via 
the Riemann-Roch theorem. Moreover, the above 
sequence also implies that 
\[
\ch_1(F) = \iota_Z^*\ch_1(E)
\]
as a sheaf on $Z$. Using a Jordan H\"older filtration and Lemma \ref{FMpoint},
it is straightforward to check that
$F$ is $\hPhi_Z-WIT_1$ and $\hPhi^1_Z(F)$ is a zero dimensional 
sheaf of length 
\[
\chi( \hPhi^1_Z(F) )= \sum_{i=1}^k r_in_i.
\]
This holds for any very ample divisor $H$ in $B$ satisfying the   genericity conditions in Definition \ref{genstabdef}. Then Lemma \ref{genvertical} implies that the set theoretical support of $\hPhi^1(E)$ is at most one dimensional. If it had dimension two or higher, the restriction of $\hPhi^1(E)$ to a generic ${\widehat Z}$ would be supported in dimension
at least one since any effective divisor on 
$\hX$ intersects ${\widehat Z}$ along 
a nonempty curve. 

Let $T \subset \hPhi^1(E)$ be the maximal zero dimensional subsheaf, 
and let $\hPhi^1(E)' = \hPhi^1(E)/T$, which is a sheaf 
of pure dimension one.  Obviously, the set theoretic support of 
$\hPhi^1(E)'$ intersects 
${\widehat Z}$ at finitely many closed points, hence
${\mathcal Tor}_1^\hX(\CO_\hZ, \hPhi^1(E)')$ is a zero dimensional 
sheaf. 
Then, using the locally free resolution
\[ 
0\to \CO_{\hX} (-{\widehat Z}) \to \CO_{\hX} \to \CO_{\widehat Z}\to 0,
\]
it follows that ${\mathcal Tor}_1^\hX(\CO_\hZ, \hPhi^1(E)')=0$
since $ \CO_{\hX} (-{\widehat Z})\otimes_\hX  \hPhi^1(E)'$ is pure of dimension one. As the higher local tor sheaves 
are obviously zero, one obtains a quasi-isomorphism 
\[ 
L\iota_\hZ^* \hPhi^1(E)'\simeq \CO_{\widehat Z}\otimes_\hX \hPhi^1(E)'.
\]
Moreover, since $T$ depends only on $E$ one can choose 
$H$ sufficiently generic such that $\hZ$ does not intersect the support 
of $T$. Then Lemma \ref{genvertical} yields an isomorphism 
\[
\CO_{\widehat Z}\otimes_\hX \hPhi^1(E)' \simeq 
\CO_{\widehat Z}\otimes_\hX \hPhi^1(E)\simeq 
\hPhi^1_Z(F).
\] 
Since ${\mathcal Tor}_i^\hX(\CO_\hZ, \hPhi^1(E)')=0$ 
for all $i\geq 1$, the Riemann-Roch theorem yields 
\[
\chi(\CO_{\widehat Z}\otimes_\hX \hPhi^1(E)') = 
\ch_2( \hPhi^1(E)')\cdot {\widehat Z}. 
\]
Therefore for any ${H}$ satisfying the required genericity assumptions there is an identity
\[
\ch_2(\hPhi^1(E))\cdot {\widehat Z}  = \chi(\hPhi_Z^1(F)) =\sum_{i=1}^k n_i r_i.  
\]
However equations \eqref{eq:FMchernB}  imply that  
\[
\ch_2(\hPhi(E)) = - \sum_{i=1}^k r_i {\hat \sigma}_*(C_i) -nf 
\]
where $\ch_3(E) = -n \ch_3(\CO_x)$. Therefore 
\[
\ch_2(\hPhi^1(E)) \cdot {\widehat Z} - \ch_2(\hPhi^0(E))\cdot {\widehat Z} = \sum_{i=1}^k n_i r_i. 
\]
In conclusion
\[ 
\ch_2(\hPhi^0(E))\cdot {\widehat Z} =0 
\]
for any very ample class $H$ in $B$. This implies that $\ch_2(\hPhi^0(E))\in \IQ \langle f\rangle $. However, equations \eqref{eq:FMchernB} imply that 
\[
\ch_i(\hPhi^0(E))=0 
\]
for $i\in \{0,1\}$ since $\hPhi^1(E)$ has one dimensional support.
Therefore $\hPhi^0(E)$ is set theoretically supported on a finite union 
of elliptic fibers. 

Now  recall that $\hPhi^0(E)$ is $\Phi-WIT_1$ 
and there is an injective morphism 
\[
\Phi^1(\hPhi^0(E)) \hookrightarrow E 
\]
according to Lemma \ref{onezero}.
Since $\hPhi^0(E)$ is $\Phi-WIT_1$ and set theoretically supported 
on a finite union of fibers, equations \eqref{eq:FMchernB} imply that $\Phi^1(\hPhi^0(E))$ will be also supported 
on a finite union of elliptic fibers. Since $E$ is 
pure of dimension two, it follows that $\Phi^1(\hPhi^0(E)) =0$, 
which further implies that $\hPhi^0(E)=0$ since $\hPhi^0(E)$ is $\Phi-WIT_1$. This implies that for any open subset of 
$\hX$ the complex $W_\bullet$ constructed in Lemma \ref{Reslemma} is a locally free resolution of $\hPhi^1(E)$. Therefore 
$\hPhi^1(E)$ must be a pure dimension one sheaf on $\hX$. 

\hfill $\Box$

For the remaining part of this section set
\be\label{eq:cpxkahlerclasses}
\omega = t\Theta-s p^* K_B, \qquad 
{\hat \omega} = {\widehat \Theta}- s {\hat p}^*K_B ,
\ee
where $s, t\in \IR$, $s> t >0$ and $s > 1$. 

Let $E$ be a vertical $\omega$-semistable sheaf on $X$
with numerical invariants 
\be\label{eq:Enuminv}
\ch_1(E)= p^*C, \qquad \ch_2(E) = kf, \qquad \ch_3(E)=-n\ch_3(\CO_x) 
\ee
where $0\neq C\in H_2(B,\IZ)$ is an effective curve class, 
$k\in (1/2)\IZ$, $n\in \IZ$. 
Suppose $E$ is also generically semistable.  Then $E$ is 
$\hPhi-WIT_1$ according to Lemma \ref{FMpuredimtwo} 
and $\hF=\hPhi^1(E)$ is a pure dimension one sheaf 
on $\hX$ with numerical invariants 
\be\label{eq:hFnuminv}
\ch_2(\hF) = {\hat \sigma}_*(C)+ n{\hat f} , \qquad 
\chi(\hF) = k - {K_B \cdot C\over 2}. 
\ee
\begin{rema}\label{positiven}
Note that Lemma \ref{FMpuredimtwo}
and  Corollary \ref{basicXcoro}.ii. imply that $n\geq 0$ for 
any sheaf $E$ as above since $\ch_2(\hF)$ must be effective. 
\end{rema} 

The next goal  is to show that 
$\hF$ is ${\hat \omega}$-semistable for sufficiently large $s$
provided that $\chi=k - K_B\cdot C /2 \geq 1$.
This will be carried out in several steps. 
For fixed 
$C,k,n$  as above with $C\neq 0$ effective, 
$\chi \geq 1$, $n\geq 0$, let  
\[
\bal 
\CS(C,k,n) =\{ & (C',l,m) \in {\rm Pic}(B)\times \IZ\times \IZ\, |\, 
  C', C-C'\, {\rm effective},\\ 
& l\geq0, \ |K_B\cdot C|l - |K_B\cdot C'|\chi\leq 0,\  0\leq m \leq n\}.
\eal 
\]
Note that $|K_B\cdot C'| \leq |K_B \cdot C|$
for any $(C',l,m) \in \CS(C,k,n)$, hence the second
 defining inequality of $\CS(C,k,n)$ yields
\[ 
0\leq l \leq \chi.
\]
Therefore $\CS(C,k,n)$ is a finite set. Moreover, 
\be\label{eq:lmineq}
|nl-m\chi|\leq n\chi
\ee
for any $(C',l,m)\in \CS(C,k,n)$. 

\begin{lemm}\label{FMstabO} 
Suppose $E$, $\hF$ are as above. Let $\hG\subset \hF$ 
be a nonzero subsheaf such that $\hF/\hG$ is a 
nonzero pure dimension one sheaf on $\hX$.
Let 
\be\label{eq:hGnuminv}
\ch_2(\hG) = {\hat \sigma}_*(C_\hG) + m {\hat  f}
\ee
with $m\in \IZ$. 
Suppose $\hG$ is ${\hat \omega}$-semistable and 
$\mu_{\hat \omega} (\hG) > 
\mu_{\hat \omega} (\hF)$ for some $s>1$. 
Then $(C_\hG, \chi(\hG), m) \in \CS(C,k,n)$. 
\end{lemm}

{\it Proof.} Given $E,\hF, \hG$ as in Lemma \ref{FMstabO}, 
note that $\mu_{\hat \omega}(\hG)>\mu_{\hat \omega}(\hF) >0$. 
Since $\hG$ is assumed ${\hat \omega}$-semistable for some 
$s>1$, Lemma 
\ref{WITO} implies that $\hG$ is $\Phi-WIT_0$. Since $E$ is $\hPhi-WIT_1$ and $\hF=\hPhi^1(E)$, Lemma \ref{WITO} implies that
$\Phi^0(\hG)$ is a subsheaf of $E$. Moreover equations \eqref{eq:FMchernA} yield  
\[
\ch_1(\Phi^0(\hG)) = p^*C_\hG, \qquad \ch_2(\Phi^0(\hG)) = 
(\chi(\hG) + K_B\cdot C_\hG/2)f, \qquad 
\ch_3(\Phi^0(\hG)) = -m\ch_3(\CO_x). 
\]
Since $\hG$ is $\Phi-WIT_0$ and $\Phi^0(\hG)$ is a nonzero 
subsheaf of $E$ one must have $C_\hG \neq 0$. Otherwise 
$\Phi^0(\hG)$ would be a nonzero sheaf 
supported on a finite union of 
elliptic fibers, leading to a contradiction since  $E$ is purely two dimensional.  
Moreover, Corollary \ref{basicXcoro} implies that $C_\hG$ 
is effective and $m\geq 0$. 
Since $\ch_2(\hF/\hG)$ must be an effective curve class,  Corollary \ref{basicXcoro}
also implies that $C = C_\hG + C'$ where
$C'$ is an effective curve class on $B$ and $n-m\geq 0$. 

Since $E$ is $\omega$-semistable, one also has $\mu_\omega(\Phi^0(\hG)) \leq \mu_\omega(E)$, 
which is equivalent to
\[
\chi(\hG)  |K_B\cdot C| - \chi(\hF)|K_B\cdot C_\hG|\leq 0. 
\]
In conclusion, $(C_\hG,\chi(\hG),m)\in \CS(C,k,n)$. 

\hfill $\Box$

Now consider the subset 
\[
\CS'(C,k,n) = \{ (C',l,m)\in \CS(C,k,n)\, |\, 
|K_B\cdot C|l - |K_B\cdot C'|\chi \leq -1\}\subset \CS(C,k,n). 
\] 
For any $s\in \IR$, $s>0$, let 
$f_s: \CS'(C,k,n) \to \IR$ be the function 
\[
f_s(C',l,m) = (s-1)(|K_B\cdot C|l - |K_B\cdot C'|\chi) + 
(nl-m\chi). 
\]
Then the following is a straightforward consequence of 
 inequality \eqref{eq:lmineq}. 
\begin{lemm}\label{sbound} 
For fixed $C,k,n$ as above there exists $s_1\in \IR$, $s_1>1$ 
depending only on  $(C,k,n)$ such that for any $s>s_1$ one has 
$f_{s}(C',l,m) <0$ for all $(C',l,m)\in \CS'(C,k,n)$. 
\end{lemm} 

\begin{lemm}\label{FMstabA}
For fixed $(C,k,n)$ as above let $s_1>1$ be as in Lemma \ref{sbound}. Then for any $s>s_1$ the Fourier-Mukai transform 
$\hF= \Phi^1(E)$ of any $\omega$-semistable and generically semistable sheaf $E$ with numerical invariants 
\eqref{eq:Enuminv} is ${\hat \omega}$-semistable.
\end{lemm}

{\it Proof.} Suppose $s>s_1$. 
The goal is to show that no destabilizing 
subsheaf $\hG \subset \hF$ as in Lemma \ref{FMstabO} can exist 
for any pair $(E,\hF)$.  Suppose $\hG\subset \hF$ is such a subsheaf
for some pair $(E, \hF)$.  Note that $\mu_{\hat \omega}(\hG) >
\mu_{\hat \omega}(\hF)$ is equivalent to 
\be\label{eq:deltaineq}
(s-1)\delta_1 + \delta_2>0,
\ee
where 
\[
\delta_1 = 
\chi(\hG)  |K_B\cdot C| - \chi(\hF)|K_B\cdot C_\hG|, \qquad 
\delta_2 = n\chi(\hG) - m \chi(\hF).
\]
According to Lemma \ref{FMstabO}, $(C_\hG, \chi(\hG), m)
\in \CS(C,k,n)$. In particular 
$\delta_1\leq 0$. Since $\delta_1\in \IZ$, there are two cases. 

$(i)$ $\delta_1\leq -1$.  Then according to Lemma \ref{sbound} 
\[
(s-1)\delta_1+\delta_2 = 
f_s(C_\hG,\chi(\hG),m)<0,
\]
contradicting \eqref{eq:deltaineq}.

$(ii)$ $\delta_1=0$. Solving for $\chi(\hG)$,
 $\delta_2$ reduces to 
\[
\delta_2= {\chi(\hF)\over |K_B\cdot C|} 
\big( n |K_B \cdot C_\hG| - m |K_B \cdot C|\big). 
\]
In this case $\mu_{\omega}(\Phi^0(\hG)) 
=\mu_{\omega}(E)$, hence one must have
\[ 
\nu_{\omega}(\Phi^0(\hG)) \leq 
\nu_{\omega}(E)
\]
since $E$ is $\omega$-semistable.
This is equivalent to $\delta_2 \leq 0$, leading again to a contradiction.

\hfill $\Box$

\subsection{From sheaves on $\hX$ to sheaves on $X$}

Again, consider K\"ahler classes of the form 
\eqref{eq:cpxkahlerclasses} on $X$, $\hX$ respectively.
Suppose $\hF$ is a pure dimension one sheaf on $\hX$
and let $L$ be a very ample line bundle on $B$. 
Using the same notation as in Definition \ref{genstabdef}, 
let $S_{\sf sm}\subset |L|$ be the nonempty open subset parametrizing 
smooth irreducible divisors $H\in |L|$ such that $Z=p^{-1}(H)$ is also 
smooth. 
Since $\hF$ is scheme theoretically supported on a closed 
subscheme of $\hX$ of pure dimension one, 
there exists a nonempty open subset $S_{\hF} \subset S_{\sf sm}$ 
such that:

$\bullet$ the set theoretic intersection between $Z_s=p^{-1}(H_s)$ 
and the support of 
$\hF$ consists of finitely many closed points, and 

$\bullet$ $H_s$ intersects the discriminant $\Delta \subset B$ transversely at finitely many points in the smooth locus of $\Delta$.

\noindent
for any closed point $s\in S_{\hF}$. If the above conditions are satisfied, 
$\hZ_s = {\hat p}^{-1}(H_s)$ is a smooth Weierstrass model 
over $H_s$ with at most finitely many nodal fibers.

\begin{lemm}\label{dimoneA}
Let $\hF$ be an ${\hat \omega}$-semistable pure 
dimension one sheaf on $\hX$ with 
\be\label{eq:positiveslope}
\ch_2(\hF) = {\hat \sigma}(C) + n{\hat f}, \qquad {\chi(\hF)}  >0,
\ee
where $C\neq 0$. 
Then the following hold. 

$(i)$ $\hF$ is $\Phi-WIT_0$ and $\Phi^0(\hF)$ is a vertical  pure dimension 
two sheaf on $\hX$. 

$(ii)$ Let $L$ be a very ample line bundle on $B$, let $H$ be a divisor 
in $B$ corresponding to a closed point in $S_{\hF}\subset |L|$, and 
 $Z=p^{-1}(H)\subset X$. Then $\Phi^0(\hF)\otimes_X \CO_Z$ is a semistable  pure dimension 
one sheaf 
on $Z$. 
 \end{lemm}

{\it Proof.} 

$(i)$ Since $\hF$ is ${\hat \omega}$-semistable,
condition \eqref{eq:positiveslope} implies that ${\rm Hom}_{\hX}(\hF, 
\CQ_x)=0$ for any closed point $x\in X$. Therefore Lemma \ref{WITO} 
implies that $\hF$ is $\Phi-WIT_0$. Moreover, equations
\eqref{eq:FMchernA} imply that
$\Phi^0(\hF)$ is a vertical  two dimensional 
sheaf. The proof of purity is completely analogous to the proof 
of Lemma \ref{Reslemma}.i.

$(ii)$ 
Under the current assumptions $\hF|_Z=\hF\otimes_\hX\CO_{\hZ}$ is a zero 
dimensional sheaf on $\hZ$. Using the same notation as in
 Lemma \ref{genvertical}, let $\Phi_{\hZ}: D^b(\hZ) \to D^b(Z)$ denote the Fourier-Mukai functor with kernel ${\cal P}|_{\hZ\times_H Z}$. 
Then it is straightforward to show that $\Phi^0_\hZ(\hF\otimes_\hX \CO_\hZ)$ is a semistable 
sheaf on $Z$ of pure dimension one set theoretically supported on a 
finite union of elliptic fibers. 
Moreover, by analogy with Lemma \ref{genvertical}, 
there is an isomorphism
\[
\Phi^0_Z(\hF\otimes_\hX\CO_{\hZ}) \simeq \Phi^0(\hF) \otimes_X 
\CO_Z 
\]

\hfill $\Box$. 

Next let $\hF$ be an ${\hat \omega}$-semistable 
pure dimension one sheaf on $\hX$ as in  Lemma \ref{dimoneA}
with 
\be\label{eq:Fnuminv}
\ch_2(\hF)= {\hat \sigma}_*(C)+ n{\hat f}, \qquad 
\chi(\hF) = k - {K_B\cdot C\over 2}\geq 1,
\ee
where $C$ is a nonzero divisor class on $B$ and $n\in \IZ$, $k\in (1/2)\IZ$. 
Note that Corollary \ref{basicXcoro} implies that $C$ must be 
effective and $n\geq 0$. 
According to Lemma \ref{dimoneA},  $\hF$ is $\Phi-WIT_0$ and $E= \Phi^0(\hF)$ 
is a vertical pure dimension two sheaf on $X$ with numerical 
invariants 
\[ 
\ch_1(E)= p^*C, \qquad \ch_2(E) = kf, \qquad \ch_3(E)=-n\ch_3(\CO_x) 
\] 
where $\CO_x$ is the structure sheaf of an arbitrary closed point $x\in X$.
In the remaining part of this section it will be shown  that $E$ is  $\omega$-semistable  for sufficiently small $t>0$. 
This will be carried out in several stages.

First suppose  $E=\Phi^0(\hF)\twoheadrightarrow G$ is a nonzero
pure dimension two quotient such that $\mu_{\omega}(G) 
\leq \mu_{\omega}(E)$ and $G$ is not isomorphic to $E$. 
Then $G$ will have numerical invariants 
\[
\ch_1(G) =p^*C_G, \qquad \ch_2(G) = \sigma_*(\alpha_G)+ cf, \qquad 
\ch_3(G) =-m\ch_3(\CO_x),
\]
where $C_G$ is a nonzero effective divisor  class on $B$, $\alpha_G$ 
is an arbitrary divisor class on $B$, and $c,m\in (1/2) \IZ$, 
$c\equiv K_B\cdot C_G/2$ mod $\IZ$,  $m\equiv C_G\cdot \alpha_G/2$ mod $\IZ$. Since $G$ is a quotient of $E$, 
not isomorphic to $E$, the curve class $C-C_G$ is effective, nonzero. 
Therefore 
\be\label{eq:CGineq} 
|K_B\cdot C_G| < |K_B \cdot C|.
\ee

\begin{lemm}\label{GboundA} 
Under the above assumptions $\alpha_G$ is an effective divisor 
class on $B$ and 
\be\label{eq:alphacbound}
{ c-|K_B \cdot \alpha_G |\over |K_B \cdot C_G | }
\leq  {k \over |K_B\cdot C|}. 
\ee
Moreover equality holds in \eqref{eq:alphacbound} if and only if 
$\alpha_G=0$. 
\end{lemm} 

{\it Proof.} Note that 
\[
\mu_{\omega}(E) = 
{1\over (s-t/2)} {k \over 
|K_B\cdot C|}, \qquad 
\mu_{\omega}(G) = {-(s/t-1)(\alpha_G \cdot K_B) + c \over 
(s-t/2)|K_B \cdot C_G|}.
\]
Given any very ample linear system $\Pi$ on $B$,
Lemma \ref{dimoneA} shows that $E|_Z$ is semistable 
for any sufficiently generic very ample divisor 
$H \in \Pi$, where $Z=p^{-1}(H)$. Moreover using Lemma 
\ref{zerotor}, one has 
\[
\chi(E|_Z)=0,\qquad \chi(G|_Z) = H \cdot \alpha_G.
\]
Therefore $H\cdot \alpha_G  \geq 0$ for any very ample 
divisor $H$ on $B$. This implies that 
$ \alpha_G$ must be an effective divisor class on $B$, 
in particular $\alpha_G \cdot K_B \leq 0$. Then 
\[
\mu_{\omega}(G) = {(s/t-1)|K_B \cdot \alpha_G| + c \over 
(s-t/2)|K_B \cdot C_G|}, 
\]
and inequality \eqref{eq:alphacbound} follows from the 
slope inequality $\mu_{\omega}(G) 
\leq \mu_{\omega}(E)$.

\hfill $\Box$

\begin{lemm}\label{GboundB} 
There exists a constant $A$ depending on $(C,k,n)$ and $s$, but not  $t$, such that 
\[ 
|c-|K_B \cdot \alpha_G || < A
\]
for all quotients $E=\Phi^0(\hF)\twoheadrightarrow G$ as above and for 
all ${\hat \omega}$-semistable sheaves $\hF$ with numerical invariants \eqref{eq:Fnuminv}. 
\end{lemm}

{\it Proof.} 
Recall that the set of isomorphism classes of ${\hat \omega}$-semistable 
sheaves with fixed numerical invariants is bounded \cite[Thm. 3.3.7]{huylehn}. Since the Fourier-Mukai transform preserves families of sheaves \cite[Prop. 6.13.]{FM_thesis}, this 
 implies that the family of 
sheaves $E=\Phi^0(\hF)$ is also bounded and depends on $(C,k,n)$, and $s$, but not $t$. 
Moreover, \cite[Lemma 1.7.6]{huylehn} implies that the same holds for the family $E_B=\sigma^*E$. 

Let $\eta_0=-K_B$, which is very ample on $B$.
Then the set of Hilbert polynomials $\calP=\{P_{\eta_0,E_B}\}$ is finite 
and indepedent of $t$. 
Let $P\in \calP$ be fixed. Obviously, the set of isomorphism 
classes $\{[E_B]\}_P$ of sheaves $E_B$ with fixed $P_{\eta_0,E_B}=P$ 
is also bounded and independent of $t$. 

Given a quotient $E\twoheadrightarrow G$,  note that 
$G_B = \sigma^*G$ is also a quotient of $E_B$, 
and there is an exact sequence of $\CO_B$-modules 
\[
0\to T_G \to G_B \to G'_B \to 0
\]
where $T_G$ is the maximal zero dimensional subsheaf of 
$G_B$ and $G'_B$ has pure dimension one.  
Since $G$ 
is pure of dimension two and has vertical support Lemma \ref{zerotor} yields an exact 
sequence 
\[
0 \to G(-\Theta) \to G \to \sigma_*G_B \to 0. 
\]
Using the above exact sequence and the Grothedieck-Riemann-Roch theorem for the embedding $\Theta \hookrightarrow X$, one obtains 
\be\label{eq:GBchernclasses}
\ch_1(G_B) = C_G, \qquad 
\ch_2(G_B) = (c-|K_B \cdot \alpha_G|)\, \ch_2(\CO_b)
\ee
with $b\in B$ an arbitrary closed point. 
Since $T_G$ is zero dimensional, $\mu_{\eta_0}(G'_B)\leq \mu_{\eta_0}(G_B)$. Then inequality \eqref{eq:alphacbound} yields 
\be\label{eq:thetaslopeA}
\mu_{\eta_0}(G'_B) \leq  {c- |K_B\cdot \alpha_G |\over 
|K_B\cdot C_G |}  \leq  {k\over | K_B\cdot C|}. 
\ee 

For fixed $P=P_{\eta_0, E_B}\in \calP$, let $\CQ_P$ denote the set of isomorphism classes of pure dimension one sheaves $F$ on $B$ such that 

$(a)$ there exists an epimorphism $E_B \twoheadrightarrow F$, 
for some $E=\Phi^0(\hF)$ as above
with $P_{\eta_0, E_B} = P$, and 

$(b)$  $\mu_{\eta_0}(F) \leq {k/ |K_B\cdot C|}$. 

Then Grothendieck's lemma 
\cite[Lemma 1.7.9]{huylehn} implies that $\CQ_P$ 
is bounded and depends only  on $P$ and the bounded family $\{[E_B]\}_P$. In particular it is independent of $t$. 
This implies that the set $\{P_{\eta_0,G'_B}\}_P$ of Hilbert polynomials of all quotients 
$E_B \twoheadrightarrow G'_B$ where $P_{\eta_0,E_B}=P$ 
is finite and $|\{P_{\eta_0,G'_B}\}_P |$ is bounded above by a constant depending 
on $P$ and  the bounded family $\{[E_B]\}_P$, but not on $t$.
Since the whole 
family $\{[E_B]\} = \cup_P \{[E_B]\}_P$ is bounded and depends only 
on $(C,k,n)$ and $s$, it follows that there exists a constant $A_1$ 
depending
on $(C,k,n)$ and $s$, but not $t$, such that 
\[
| \chi(G'_B)| < A_1 
\]
for all pure dimension one quotients $E_B \twoheadrightarrow G'_B$, 
for all $E=\Phi^0(\hF)$ as above.  

To conclude the proof, note that 
$\chi(G_B) = \chi(T_G) + \chi(G'_B) \geq \chi(G'_B)$ 
since $T_G$ is zero dimensional. On the other hand, using
equation \eqref{eq:GBchernclasses} and the Riemann-Roch theorem,
\[
\chi(G_B) = c-|K_B\cdot \alpha_G|+|K_B\cdot C_G|/2.
\]
Therefore, using inequality \eqref{eq:CGineq}, 
\[
 c-|K_B\cdot \alpha_G | > -A_1 - |K_B\cdot C_G|/2> 
-A_1 - |K_B\cdot C|/2.
\]
At the same time inequalities \eqref{eq:CGineq},  \eqref{eq:alphacbound}
yield
\[ 
c-|K_B\cdot \alpha_G | \leq {|K_B\cdot C_G|\over |K_B\cdot C|}k<|k|. 
\]
Therefore the claim follows.

\hfill $\Box$

\begin{lemm}\label{FMstabB}
There exists a constant $t_1\in \IR$, $0<t_1<s$,
 depending on $(C,k,n)$ and $s$,
such that for all $0<t<t_1$ and for any ${\hat \omega}$-semistable 
 sheaf  $\hF$ with numerical invariants \eqref{eq:Fnuminv}
any pure dimension two quotient $E\twoheadrightarrow G$ with $\mu_{\omega}(G) \leq \mu_{\omega}(E)$ 
is vertical, where $E=\Phi^0(\hF)$.
\end{lemm} 

{\it Proof.} For any $0<t<s$ set 
$\omega_t = t \Theta -sp^*K_B$.
The numerical invariants $(C,k,n)$ and $s>0$ are fixed in the following. 

Suppose the opposite statement holds. 
Given any $0<t_1<s$, there exist $0<t<t_1$, 
a sheaf $\hF$ as in Lemma \ref{FMstabB} 
and a nonzero quotient $E\twoheadrightarrow G$, 
not isomorphic to $E$, such that $\mu_{\omega_t}(G) \leq \mu_{\omega_t}(E)$ and $G$ is not vertical. It will be shown below that this leads to a contradiction.  

Note that $G$ has numerical invariants 
\[
\ch_1(G) =p^*C_G, \qquad \ch_2(G) = \sigma_*(\alpha_G)+ cf, \qquad 
\ch_3(G) =-m\ch_3(\CO_x)
\]
and $G$ is vertical if and only if $\alpha_G =0$. Suppose $\alpha_G\neq 0$. 
Lemma \ref{GboundA} shows that $\alpha_G$ is effective, hence 
\be\label{eq:GEslopes}
\mu_{\omega}(G) = {s\over t(s-t/2)} 
{{|K_B\cdot \alpha_G }|\over {|K_B\cdot C_G |} }
+ \delta, \qquad 
\delta = {1\over s-t/2} {c-|K_B \cdot \alpha_G|\over |K_B\cdot C_G 
|}.
\ee
According to Lemma \ref{GboundB}, there is a 
constant $A$ depending on $(C,k,n)$ and  $s$, but not $t$, 
such that 
\[ 
|c -|K_B \cdot \alpha_G|| <A
\] 
for any  quotient $E\twoheadrightarrow G$ as above. 
Moreover since $-K_B$ is very ample, the set 
\[
\{|\beta \cdot K_B|\, |\, 
0\neq \beta \in {\rm Pic}(B) \ {\rm effective}\} \subset \IZ_{>0}
\]
is bounded from below. Let $M\in \IZ_{>0}$ denote its minimum
and note that 
$| K_B\cdot C_G|\geq M$, $|K_B \cdot \alpha_G|\geq M$ since $C_G, \alpha_G$ are effective, nonzero. 

Suppose $0<t<2$. Then 
$0<s-1< s-t/2$, hence 
\[
|\delta| < {1\over s-1} {A\over M}.
\]
Using inequality \eqref{eq:CGineq},
\[
{s\over t(s-t/2)} {
{|K_B\cdot \alpha_G }|\over {|K_B\cdot C_G |} } 
> {s\over t(s-t/2)} {M\over |K_B \cdot C|}.
\]
Moreover, the map 
\[ 
f: (0,\ s)\to \IR, \qquad f(t)= {s\over t(s-t/2)} 
\]
is a decreasing function of $t$ on the interval $0<t<s$ for fixed $s>0$,
and $\lim_{t\to 0}  f(t)=+\infty$. 
Therefore there exists a constant $0<t_1< \rm{min}\{s,2\}$
depending on $(C,k,n)$ and $s$ such that for any $0<t<t_1$, 
\[
\mu_{\omega}(G) > {1\over (s-1)} {|k|  \over 
|K_B\cdot C|} +1 
\]
for all quotients $G$ as in Lemma \ref{FMstabB} with 
$\alpha_G\neq 0$. In order to conclude the proof note that 
under the current assumptions 
\[ 
|\mu_{\omega}(E)| = {|k|\over 
(s-t/2) |K_B\cdot C|} 
\leq {1\over (s-1)} {|k| \over 
|K_B\cdot C|}
\]
for any $0<t<t_1$, leading to a contradiction. 

\hfill $\Box$

\begin{lemm}\label{FMstabC}
Let $s>s_1$ be fixed, where $s_1>1$ is a constant as in Lemma \ref{sbound} and $0<t <t_1$ where $t_1$ is a constant as
in Lemma \ref{FMstabB} for fixed $(C,k,n)$ and $s$. 
Then the Fourier-Mukai transform $E=\Phi^0(\hF)$ of any 
${\hat \omega}$-semistable 
 sheaf  $\hF$ with numerical invariants \eqref{eq:Fnuminv} is 
$\omega$-semistable for all $0<t<t_1$. 
\end{lemm} 

{\it Proof.} Recall that under the current assumptions 
\[
\chi(\hF)=k-{K_B\cdot C\over 2}\geq 1.
\]
Let $E \twoheadrightarrow G$ be a nonzero pure dimension two 
quotient of $E$ such that $G$ is $\omega$-semistable 
and destabilizes $E$. This means either
\[
\mu_{\omega}(G) < \mu_{\omega}(E) = 
{1\over (s-t/2)} {k \over 
|K_B\cdot C|}
\]
or $\mu_{\omega}(G) = \mu_{\omega}(E)$ 
and 
\[ 
 \nu_{\omega}(G) < \nu_{\omega}(E). 
\]
According to Lemma \ref{FMstabB}, $G$ must be 
vertical i.e. 
\[
\ch_1(G) =p^*C_G, \qquad \ch_2(G) = cf, \qquad 
\ch_3(G) =-m\ch_3(\CO_x)
\]
where $C_G$ is a nonzero effective divisor class on $B$
 and $c\in (1/2) \IZ$, $m\in \IZ$, 
$c\equiv K_B\cdot C_G/2$ mod $\IZ$. 
Therefore
\[ 
\mu_{\omega}(G) = {c \over 
(s-t/2)|K_B\cdot C_G|}.
\]
At the same time 
$E$ is generically semistable according to Lemma \ref{dimoneA}. Hence, given any very ample linear system $\Pi$ on $B$, 
$E|_Z$ is semistable 
for any sufficiently generic very ample divisor 
$H \in \Pi$, where $Z=p^{-1}(H)$. Moreover Lemma \ref{zerotor} 
yields $\chi(E|_Z)=
\chi(G|_Z)=0$. This implies that $G$ must be 
generically semistable as well. Then Lemma \ref{FMpuredimtwo} implies that $G$ is $\hPhi-WIT_1$ and $\hPhi^1(G)$ is pure dimension one. Furthermore the
epimorphism $E\twoheadrightarrow G$ yields an epimorphism 
$\hF \twoheadrightarrow \hPhi^1(G)$. Therefore 
\be\label{eq:hatslopes}
\mu_{{\hat \omega}}(\hF) \leq 
\mu_{{\hat \omega}}(\hPhi^1(G))
\ee
since $\hF$ is ${\hat \omega}$-semistable. 
The numerical invariants of $\hPhi^1(G)$ are 
\[
\ch_2(\hPhi^1(G)) = {\hat \sigma}_*(C_G) + m{\hat f}, \qquad 
\chi(\hPhi^1(G)) = 
c - K_B\cdot C_G  /2.
\]
Note that Corollary \ref{basicXcoro} implies that $0\leq m \leq n$.
Moreover, $\chi(\hPhi^1(G))>0$ since $\chi(\hF)>0$
under the current assumptions. 
At the same time, the slope inequality $\mu_{\omega}(G)\leq 
\mu_{\omega}(E)$ is equivalent to $\delta_1 \leq 0$,
where 
\[
\delta_1 = \chi(\hPhi^1(G)) |K_B\cdot C| - \chi(\hF) 
|K_B \cdot C_G |.
\]
In conclusion $(C_G, \chi(\hPhi^1(G)), m) \in \CS(C,k,n)$, 
where $\CS(C,k,n)$ is the finite set defined above Lemma 
\ref{FMstabO}. 

The slope inequality \eqref{eq:hatslopes} is equivalent to
\[
(s-1) \delta_1 + \delta_2 \geq 0,
\]
where 
\[
\delta_2 = n \chi(\hPhi^1(G)) 
- m\chi(\hF).
\]

Since $\delta_1\in \IZ$, one has to distinguish two cases.  

$(i)$ $\delta_1\leq -1$. In this case 
\[
(s-1)\delta_1+\delta_2 = f_s(C_G, \chi(\hPhi^1(G)), m) <0 
\]
where $f_s:\CS'(C,k,n)\to \IR$ is the function defined above Lemma \ref{sbound}. 
Obviously, this leads to a contradiction. 

$(ii)$ Suppose $\delta_1=0$. This implies 
\[
\delta_2=  {\chi(\hF)\over |K_B\cdot C|}
\big( n |K_B \cdot C_G| - m |K_B \cdot C|\big). 
\]
However in this case $\mu_{\omega}(G) 
=\mu_{\omega}(E)$, hence one must have
\[ 
\nu_{\omega}(G) < 
\nu_{\omega}(E), 
\]
which is equivalent to $n |K_B \cdot C_G| - m |K_B \cdot C|<0$. 
Since $\chi(\hF)>0$, this implies $\delta_2 < 0$, leading again to 
 a contradiction. 

\hfill $\Box$

\subsection{Proof of Theorem \ref{main}}

This subsection concludes the proof of Theorem \ref{main}. 
Let ${\hat \gamma} \in H_2(B,\IZ) \oplus \IZ\oplus \IZ$ be
fixed numerical invariants with ${\hat \gamma}_1$ an 
effective curve class on $B$, and ${\hat \gamma}_3>0$. 
Let $\gamma = \phi({\hat \gamma})$. Let $s_1({\hat \gamma})>1$ be a constant as 
in Lemma \ref{sbound}. For any $s\in \IR$, $s>s_1({\hat \gamma})$, let $t_1({\hat \gamma},s) \in \IR$, $0< t_1({\hat \gamma},s) <s$ be a
constant as in Lemma \ref{FMstabB}. Let 
\[
\omega = t\Theta - sp^*K_B, \qquad {\hat \omega} = {\widehat \Theta}
- s{\hat p}^*K_B. 
\]

Lemmas \ref{dimoneA} and \ref{FMstabC} prove that any ${\hat \omega}$-semistable sheaf $\hF$ with numerical invariants ${\hat \gamma}$ is $\Phi-WIT_0$ and $\Phi^0(\hF)$ is an $\omega$-semistable vertical pure dimension two sheaf $E$ on $X$ with invariants 
$\gamma$. Moreover $E$ is 
also generically semistable. 
Conversely, Lemmas \ref{FMpuredimtwo} and \ref{FMstabA} 
prove that any $\omega$-semistable and generically semistable 
vertical sheaf $E$ with numerical invariants $\gamma$ is 
${\hPhi}-WIT_1$ and $\hPhi^1(E)$ is an ${\hat \omega}$-semistable 
sheaf on $\hX$ with invariants ${\hat \gamma}$. 
Furthermore, 
Lemmas \ref{genstablemma} and \ref{adiabstablemma} prove that 
generic semistability is equivalent to adiabatic semistability
for $\omega$-semistable sheaves.

In order to conclude the proof of Theorem \ref{main}.i. note 
that the Fourier-Mukai transform preserves flat families of 
sheaves \cite[Prop. 6.13.]{FM_thesis}. 

For the second statement, note that the substack $\CM^{\sf ad}_\omega(X,\gamma)$ is open in $\CM_\omega(X,\gamma)$ according to Lemma 
\ref{adiabopen}. Moreover, let $M_{\hat \omega}(\hX,{\hat\gamma})$, 
$M_\omega(X,\gamma)$ be the coarse moduli schemes parameterizing 
$S$-equivalence classes of semistable sheaves. 
As noted above Lemma \ref{adiabopen}, according 
to \cite[Ex. 8.7]{goodmoduli}, the coarse moduli schemes are 
good moduli coarse moduli spaces for the moduli stacks
$\CM_{\hat \omega}(\hX, {\hat \gamma})$, $\CM_\omega(X,\gamma)$. 
Using \cite[Thm 4.16]{goodmoduli}, this yields a commutative diagram 
\[
\xymatrix{ 
\CM_{\hat \omega}(\hX, {\hat \gamma})\ar[r]^-{\varphi} 
\ar[d]^-{\hat \varrho} & 
\CM_\omega(X,\gamma) \ar[d]^-{\varrho} \\
M_{\hat \omega}(\hX,{\hat\gamma})\ar[r]^-{f} & M_\omega(X,\gamma)
\\ }
\]
where $\varphi$ factors through the natural embedding 
$\CM_\omega^{\sf ad}(X, \gamma) \subset \CM_\omega(X,\gamma)$. 
In the above diagram $f$ is a morphism of schemes, and the vertical morphisms 
are those constructed in  \cite[Thm 4.16]{goodmoduli}. 
Since both coarse moduli spaces are 
projective, it follows that $f$ is proper. At the same time, 
according to Lemma \ref{adiabopen}, 
the scheme theoretic image $M_\omega^{\sf ad}(X,\gamma)= \varrho(\CM_{\omega}^{\sf ad}(X,\gamma))$ 
is open in $M_\omega(X,\gamma)$. Since $f$ is proper and 
$M_{\hat \omega}(\hX, {\hat \gamma})$ is projective, it follows 
that $M_\omega^{\sf ad}(X,\gamma)$ is open and closed in 
$M_\omega(X,\gamma)$. Therefore $\CM_\omega^{\sf ad}(X,\gamma)$ is open and closed in $\CM_\omega(X,\gamma)$. 

\hfill $\Box$

\section{Vertical sheaves on elliptic $K3$ pencils}\label{verticalKthree}

Using the notation in Section \ref{countingsect}, 
let $X$ be a smooth generic Weierstrass model over the Hirzebruch 
surface $B = \IF_a$, 
$0\leq a \leq 1$. Let $\pi: X \to \IP^1$ be the natural projection to 
$\IP^1$. Note that all fibers of $\pi$ are reduced irreducible 
elliptic $K3$-surfaces 
in Weierstrass form. For sufficiently generic $X$, the generic $K3$ fiber 
is a smooth Weierstrass model and the singular fibers will be Weierstrass 
models with finitely many isolated type $I_1$ and 
$I_2$ fibers. In particular 
all singular K3 fibers are reduced, irreducible with isolated simple nodal 
singularities. 
This will be assumed throughout this section. 

Let 
$\Xi\in {\rm Pic}(B)\simeq H_2(B,\IZ)$ 
denote the fiber class of the Hirzebruch surface and note that 
the $K3$ fiber class is $D=p^*\Xi\in {\rm Pic}(X)\simeq H_4(X,\IZ)$. 
Let $\omega = t\Theta - sp^*K_B$ 
be a K\"ahler class on $X$ with $t,s\in \IR$, $0<t<s$. 
In order to simply the notation the pushforward $\sigma_*C\in H_2(X,\IZ)$ of a curve class on $B$ will be denoted by $C$. The distinction 
will be clear from the context. 

First note the following simple fact. 
\begin{lemm}\label{chernclasseslemma}
Let $E$ be a nonzero pure dimension two sheaf set theoretically supported 
on a finite union of $K3$ fibers. Then 
\be\label{eq:chernclassesA} 
\ch_1(E) = rD, \qquad \ch_2(E) = m\Xi+ lf
\ee 
for some $r,m,l\in \IZ$, $r\geq 1$.  
\end{lemm} 

{\it Proof.} Let $C_0\in H_2(B,\IZ)$ be a section class with $C_0^2 =-a$ 
for $B = \IF_a$, $0\leq a\leq 1$. Let $D_0= p^*C_0 \in H_4(X,\IZ)$. 
According to Lemma \ref{basicX} 
\[
H_4(X,\IZ)/_{\rm torsion} \simeq \IZ\langle D_0,D, \Theta\rangle,\qquad 
H_2(X,\IZ)/_{\rm torsion} \simeq \IZ\langle C_0, \Xi, f\rangle. 
\]
Obviously $\ch_i(E) \cdot D =0$ in the intersection ring of $X$
for $1\leq i\leq 3$ 
since under the current assumptions 
$E\simeq E\otimes_X \CO_X(-D)$ and $D^2=0$.
Then the claim follows easily from Lemma \ref{basicX} 
and the following 
relations in the intersection ring  of $X$:
\be\label{intersrel}
\bal 
& D_0\cdot D = f, \qquad D\cdot D =0, \qquad \Theta \cdot D = \Xi\\
& 
C_0 \cdot D =1, \qquad \Xi\cdot D =0, \qquad f \cdot D =0.\\
\eal
\ee

\hfill $\Box$ 

Now suppose $\iota: S\hookrightarrow X$ is a singular K3 fiber. Under the current genericity assumptions $S$ is an elliptic surface over $\IP^1$ 
in Weierstrass form with finitely 
many type $I_1$ and $I_2$ fibers. Therefore $S$ will have 
finitely many isolated simple nodes and the singular locus of $S$ 
is disjoint from the canonical section of the Weierstrass model. 
Let $\rho: {\widetilde S} \to S$ be a smooth 
crepant resolution of singularities and let $\phi=\iota\circ \rho : {\widetilde S} \to X$. Let ${\widetilde \Xi} = \phi^*\Theta$, ${\tilde f} = \phi^*D_0$ be 
induced divisor classes on ${\widetilde S}$. Let also 
$\epsilon_1, \ldots, \epsilon_k$ denote the exceptional 
$(-2)$ curve classes on $\wS$ and note that 
\be\label{eq:resinters}
{\widetilde \Xi}^2 = -2, \qquad {\widetilde \Xi}\cdot {\tilde f} =1, 
\qquad {\tilde f}^2 =0, \qquad 
\epsilon_i \cdot {\widetilde \Xi} = \epsilon_i \cdot {\tilde f} =0, \qquad 
1\leq i \leq k, 
\ee
in the intersection ring of ${\widetilde S}$, and 
\[ 
\phi_*\epsilon_i =0, \qquad 1\leq i \leq k.
\]
Then note the following.

\begin{lemm}\label{pushlemma} 
$(i)$ Let $\iota: S \hookrightarrow$ be a smooth $K3$ fiber of $X$
and $F$ a torsion free sheaf on $S$ such that 
\[
\ch_1(\iota_*F) = rD, \qquad \ch_2(\iota_* F) = m\Xi + l f, \qquad 
\ch_3(\iota_*F) = -n \ch_3(\CO_x) 
\]
for some $l,m,n,r\in \IZ$, $r\geq 1$. Then 
\[
\ch_0(F) =r, \qquad \ch_1(F) = m\Xi+lf + \beta, \qquad 
\ch_2(F) = -n \ch_2(\CO_s) 
\]
for a curve class $\beta \in H_2(S,\IQ)$ such that $\beta\cdot \Xi=\beta\cdot f=0$ in the intersection ring of $S$ and 
$\iota_* \beta=0$. 

$(ii)$ Let $\iota: S \hookrightarrow$ be a singular $K3$ fiber of $X$
and $\wF$ a torsion free sheaf on the resolution $\wS$ such that 
\[
\ch_1(\phi_*\wF) = rD, \qquad \ch_2(\phi_* \wF) = m\Xi + l f, \qquad 
\ch_3(\phi_*\wF) = -n \ch_3(\CO_x) 
\]
for some $l,m,n,r\in \IZ$, $r\geq 1$. Then 
\[
\ch_0(\wF) =r, \qquad \ch_1(F) = m\wXi+l{\tilde f }+ \tbeta + 
\sum_{i=1}^k p_i \epsilon_i, \qquad 
\ch_2(\wF) = -{\tilde n} \ch_2(\CO_s) 
\]
for some $p_i \in \IQ$, $1\leq i \leq k$, ${\tilde n}\in \IZ$, ${\tilde n}\leq n$, and a curve class $\tbeta \in H_2(\wS,\IQ)$ such that \[
\tbeta\cdot \wXi=\tbeta\cdot {\tilde f}=\tbeta\cdot \epsilon_i =0, 
\qquad 1\leq i \leq k
\]
in the intersection ring of $\wS$ and $\phi_*\tbeta=0$. 
\end{lemm} 

{\it Proof.} For $(i)$ note that the Grothendieck-Riemann-Roch 
theorem yields 
\be\label{eq:GRRa}
\ch_0(F) =r, \qquad \iota_*\ch_1(F) = m\Xi+lf, 
\qquad \ch_2(F) =-n \ch_3(\CO_s)
\ee
with $s\in S$ a closed point. 
Then the push pull formula yields 
\[ 
\ch_1(F) \cdot \Xi = l-2m, \qquad 
\ch_1(F)\cdot f = m
\]
in the intersection ring of $S$. Therefore 
\[ 
\ch_1(F) = m\Xi + lf + \beta 
\]
where $\beta \in H_2(S,\IZ)$ is orthogonal to $\Xi, f$. 
Moreover the second equation in \eqref{eq:GRRa} 
implies $\iota_*\beta =0$. 

$(ii)$ Since 
$\phi = \iota \circ \rho$ and $\rho: \wS\to S$ is an isomorphism onto the smooth open part of $S$, $R^1\phi_*\wF$ is a zero dimensional sheaf supported at the nodes of $S$. Then the  
Grothendieck-Riemann-Roch theorem gives 
\[
\ch_0(\wF)=r, \qquad \phi_*\ch_1(\wF) = m\Xi+l f, \qquad 
\phi_*\ch_2(\wF) = -n \ch_3(\CO_x) + \ch_3(R^1\phi_*\wF).
\]
The remaining part of the proof is analogous to $(i)$. 

\hfill $\Box$ 

For any pure dimension two sheaf $E$ with scheme theoretic support 
on a reduced nodal K3 fiber $S\subset X$, let $\wF_E = \phi^*E/{\rm torsion}$. Note that given an ample class 
 $\omega$ on $X$, the real divisor class 
${\tilde \omega}_\lambda = \lambda \phi^*(\omega|_S)-\sum_{i=1}^k \epsilon_i$ is ample
on $\wS$ for sufficiently large $\lambda \in \IR$, $\lambda >0$. Then the following result is similar to \cite[Lemma 2.1]{ModuliReflexive}.

\begin{lemm}\label{stabreflexive} 
Let $\iota:S \to X$ be a reduced nodal $K3$ fiber. 
Let $E$ be a nonzero $\omega$-slope stable 
pure dimension two sheaf  on $X$ 
set theoretically supported on $S$.
Then $E$ is  
scheme theoretically  supported  on $S$ and 
$\wF_E$ is ${\tilde \omega}_\lambda$-slope stable for sufficiently
large $\lambda>0$.  
\end{lemm} 

{\it Proof.} Proving that $E$ is scheme theoretically  supported  on $S$ 
is completely 
analogous to Lemma \ref{genlemmaA}. The details will be omitted.
For the second statement, by construction there is an exact sequence 
\[ 
0 \to \CT \to \phi^*E \to \wF_E \to 0
\]
where $\CT$ is set theoretically supported on the exceptional 
locus of $\rho$. This yields a second sequence 
\[ 
0 \to \phi_*\CT \to \phi_*\phi^*E {\buildrel f\over \longto}  \phi_*\wF_E \to 
R^1\phi_*\CT\to \cdots 
\]
where $ \phi_*\CT, R^1\phi_*\CT$ are set theoretically supported 
on the singular locus $S^{\sf sing}\subset S$, which consists of 
fintely many points. Moreover there is a natural morphism 
$g:E \to \phi_*\phi^*E$ which is an isomorphism on the smooth 
locus $S\setminus S^{\sf sing}$. 
The morphism $f\circ g : E \to \phi_*\wF_E$ is also an isomorphism 
on $S\setminus S^{\sf sing}$, hence it must be injective since $E$
is purely two dimensional.  In conclusion there is an exact sequence 
\be\label{eq:seqone} 
0\to E \to \phi_*\wF_E {\buildrel f\over \longto} T\to 0
\ee 
with $T$ zero dimensional. 
This implies that $\mu_\omega(E) = \mu_\omega(\phi_*\wF_E)$. 

If $r=1$, $\wF_E$ is a rank one torsion free sheaf which is slope stable for any polarization of $\wS$. Recall that slope stability is defined 
with respect to saturated nonzero test subsheaves as in \cite[Sect. 1.6]{huylehn}.

Let $r\geq 2$ and suppose $\wG \subset \wF_E$ is a nonzero proper saturated subsheaf
of rank $1\leq r'\leq r-1$. Then $\phi_*\wG$ is a subsheaf of $\phi_*\wF_E$. Let $I\subset T$, $G\subset \phi_*\wG$ be the image and respectively the kernel of  $f|_{\phi_*\wG}$ in the exact sequence 
\eqref{eq:seqone}. Then $I$ is zero dimensional and $G$ is a subsheaf of $E$. 
This implies that 
$\mu_\omega(G) = \mu_\omega(\phi_*\wG)$, hence 
$\mu_\omega(\phi_*\wG) < \mu_\omega(E)=\mu_\omega(\phi_*\wF_E)$ since $E$ is $\omega$-stable by assumption. 
Therefore 
\be\label{eq:diffA}
(r'\ch_1(\wF_E) -r\ch_1(\wG)) \cdot \phi^*\omega >0.
\ee

Let $\lambda_0 >0$ be fixed such that ${\tilde \omega}_0= 
\lambda_0\phi^*\omega - \sum_{i=1}^k\epsilon_i$ 
is ample on $\wS$. The subsheaves $\wG \subset \wF_E$ 
are of two types:

$a)$ $(r'\ch_1(\wF) -r\ch_1(\wG)) \cdot \omega_0>0$. 
Then, using inequality \eqref{eq:diffA}, 
\be\label{eq:diffB}
(\omega_0 + \lambda\phi^*\omega) \cdot (r'\ch_1(\wF_E) - 
r\ch_1(\wG) >0 
\ee
for any $\lambda >0$. 

$b)$ $(r'\ch_1(\wF_E) -r\ch_1(\wG)) \cdot \omega_0\leq 0$. 
According to Grothendieck's Lemma \cite[Lemma 1.7.9]{huylehn} 
the family of such subsheaves is bounded for fixed $\wF_E$ and $\omega_0$. Therefore there exists a constant $c_1>0$ 
depending on 
$\wF_E$, $\omega_0$ such that 
\[
(r'\ch_1(\wF_E) -r\ch_1(\wG)) \cdot \phi^*\omega >c_1
\]
for any subsheaf $\wG$ of type $(b)$. Furthermore there is a 
second constant $c_2>0$ depending on 
$\wF$, $\omega_0$ such that 
\[
(r'\ch_1(\wF_E) -r\ch_1(\wG)) \cdot \omega_0> -c_2
\]
for any such subsheaf. This implies that there exists a sufficiently large $\lambda>0$ such that 
inequality \eqref{eq:diffB} holds for all subsheaves of type $(b)$
as well.  
In conclusion $\wF_E$ is $(\omega_0 + \lambda\phi^*\omega)$-slope 
stable. 

\hfill $\Box$

Now recall that the discriminant of a rank $r\geq 1$ torsion 
free sheaf $F$ 
on a smooth projective surface $S$ is defined (up to normalization) 
by 
\[
\Delta(F) = n+ {1\over 2r} \ch_1(F)^2  
\]
where $\ch_2(F) = -n \ch_2(\CO_s)$, with $s\in S$ and arbitrary closed 
point. For any  vertical pure dimension two sheaf $E$ 
with 
\[
\ch_1(E)=rD, \qquad \ch_2(E)=m\Xi+lf, \qquad \ch_3(E)=-n\ch_3(\CO_x)
\]
let 
\be\label{eq:discrimA}
\delta(E) = n - {1\over r} m(m-l).
\ee
Then note the following. 

\begin{lemm}\label{bogomolovlemma}
Let $E$ be an 
$\omega$-slope semistable pure dimension two sheaf on $X$ with numerical invariants 
\[ 
\ch_1(E) = rD, \qquad \ch_2(E) = m\Xi+lf, \qquad \ch_2(E) =-n \ch_3(\CO_x)
\]
where $r,l,m,n\in \IZ$, $r\geq 1$, and $x\in X$ is an arbitrary 
closed point. Suppose $E$ is scheme theoretically supported on  a reduced $K3$ fiber 
$\iota: S \hookrightarrow X$. Then 
$\delta(E)\geq 0$. 
\end{lemm}

{\it Proof.}
Obviously, $E=\iota_*F$ for a torsion free sheaf on $S$. 

Suppose first that  $S$ is smooth. Then $F$ is $\omega|_S$-slope semistable. According to Lemma \ref{pushlemma}.i,
\[
\ch_1(F) = m \Xi + lf + \beta
\]
where $\beta \in H_2(S,\IQ)$ is a curve class 
such that $\beta\cdot \Xi = \beta \cdot f =0$. 
At the same time $\omega|_S = t \Xi + 2sf$, hence 
$\beta \cdot \omega|_S=0$. Then
$\beta^2 \leq 0$ according to the  Hodge index theorem. 
Since $F$ is $\omega|_S$-slope semistable,  it satisfies the Bogomolov inequality,
$\Delta(F) \geq 0$, where  
\[
\Delta(F) =n -{1\over r} m(m-l)  + {\beta^2 \over 2r} = \delta(E) +{\beta^2\over 2r}.
\]
Since $\beta^2\leq 0$, this implies  the claim. 

Next let $S$ be a singular $K3$ fiber. Suppose first that $E$ is $\omega$-slope stable. 
Then it is scheme theoretically supported on $S$. 
Let $\wF_E$ be the corresponding torsion free sheaf on $\wS$. 
Lemma \ref{stabreflexive} shows that $\wF_E$ is stable for a suitable ample class ${\tilde \omega}$ 
on $\wS$, hence $\Delta(\wF_E)\geq 0$. 
Moreover as shown in the proof of Lemma \ref{stabreflexive}, there is an exact sequence 
\[ 
0 \to E \to \phi_*\wF_E \to T \to 0
\] 
with $T$ zero dimensional. Setting 
\[
\ch_3(E) = -n \ch_3(\CO_x), \qquad \ch_3(\phi_*\wF_E) = 
- n' \ch_3(\CO_x)
\]
this implies $n \geq n'$. Furthermore, 
according to Lemma \ref{pushlemma}.ii, 
\[
\ch_0({\widetilde F}_E) = r,\qquad 
\ch_1(\wF_E) = m {\widetilde \Xi} + l{\tilde f}+  {\tilde \beta}+ 
\sum_{i=1}^k p_i \epsilon_i, \qquad 
\ch_2(\wF_E) = -{\widetilde n} \ch_2(\CO_s)
\]
with $p_i \in \IQ$, ${\tilde n}\in \IZ$, ${\tilde n}\leq n'$, 
and ${\tilde \beta}\in H_2({\widetilde S}, \IQ)$ 
a curve class orthogonal to ${\widetilde \Xi},  {\tilde f},
\epsilon_i$ for all $1\leq i\leq k$. 
In particular ${\tilde \beta} \cdot {\tilde \omega} =0$. 
Then  
\[
\Delta({\widetilde F}_E) = {\tilde n}- {1\over r} m(m-l) + {1\over 2r} 
\big( {\tilde \beta}^2 -2\sum_{i=1}^n p_i^2 \big) 
\]
Since ${\tilde \beta} \cdot {\tilde \omega}=0$,  the Hodge index theorem shows that 
${\tilde \beta}^2 \leq 0$. 
Since ${\tilde n}\leq n'\leq n$, this implies the claim. 

To finish the proof, suppose $E$ is strictly $\omega$-slope semistable. According to \cite[Thm 1.6.7.ii]{huylehn}, there is a 
 Jordan-H\"older filtration 
\[
0=E_0 \subset E_1 \subset \cdots \subset E_j =E
\]
for slope semistability with $j\geq 2$. Each succesive quotient 
$E_i/E_{i-1}$, $1\leq i \leq j$, is $\omega$-slope polystable, hence scheme theoretically supported on $S$. Therefore 
$\delta(E_i/E_{i-1})\geq 0$ for all $1\leq i \leq j$. 
Then the claim follows by a recursive application of Lemma 
\ref{extlemma} below. 

\hfill $\Box$

\begin{lemm}\label{extlemma} 
Let 
\[
0\to E_1\to E\to E_2\to 0
\]
be an extension of nonzero pure dimension two sheaves such that 
$E_1,E_2$ are $\omega$-slope semistable and set theoretically supported on finite unions of K3 fibers.
Suppose that $\mu_\omega(E_1)=
\mu_\omega(E_2)$. Then 
\[
\delta(E) \geq \delta(E_1)+\delta(E_2) 
\]
\end{lemm}

{\it Proof}. Let 
\[
\ch_1(E_i)=r_iD,\qquad \ch_2(E_i) = m_i \Xi+l_il
\]
for $1\leq i \leq 2$, where $r_1, r_2\geq 1$. Then 
\[
\delta(E)-\delta(E_1)-\delta(E_2) = d
\]
where 
\[
d = {(r_1m_2-r_2m_1)\over r_1r_2(r_1+r_2)} 
\left[(r_1m_2-r_2m_1)-(r_1l_2-r_2l_1)\right].
\]
Let $S$ be a generic smooth K3 fiber and $\Xi, f\in H_2(S,\IZ)$ 
the section, respectively fiber class. Then 
\[ 
d = -{\alpha^2 \over 2r_1r_2(r_1+r_2)}
\]
where 
\[
\alpha = r_1(m_2\Xi+l_2f) - r_2(m_1\Xi+l_1f).
\]
The slope equality $\mu_\omega(E_1)=\mu_\omega(E_2)$ is equivalent to $\alpha \cdot \omega|_S=0$. Since $S$ is smooth, the
Hodge index theorem shows that $\alpha^2\leq 0$. This proves the claim.

\hfill $\Box$

\begin{lemm}\label{wallbound}
Let $E$ be an 
$\omega$-slope semistable sheaf on $X$ with numerical invariants 
\[ 
\ch_1(E) = rD, \qquad \ch_2(E) = m\Xi+lf, \qquad \ch_2(E) =-n \ch_3(\CO_x)
\]
where $r,l,m,n\in \IZ$, $r\geq 1$, and $x\in X$ is an arbitrary 
closed point. Suppose $E$ is scheme theoretically supported on a 
reduced $K3$ fiber $\iota:S\hookrightarrow X$ and there is an extension
\[
0\to E_1 \to E \to E_2 \to 0
\]
with $E_1, E_2$ nonzero 
pure dimension two sheaves with $\ch_1(E_i) = r_iD$, $r_i \in \IZ$, $r_i\geq 1$, $1\leq i\leq 2$. Moreover suppose
\be\label{eq:extconditions}
\mu_\omega(E_1) = \mu_\omega(E_2) \quad {\rm and}\quad
{1\over r_1} \ch_2(E_1) - {1\over r_2} \ch_2(E_2) \neq 0.
\ee
Then 
\be\label{eq:wallboundA}
{t\over s} \geq {2\over 1+ r^3\delta(E)}.
\ee
\end{lemm} 

{\it Proof.} 
As in the proof of Lemma \ref{extlemma}, let 
 $\iota':S'\hookrightarrow X$ be a smooth generic K3 fiber and 
$\Xi,f\in H_2(S', \IZ)$ the section and fiber class respectively. 
Note that 
\[ 
\ch_2(E_i) = \iota'_* \alpha_i
\]
for $\alpha_i = m_i \Xi+l_i f\in H_2(S', \IZ)$, $1\leq i\leq 2$
and 
\[ 
\delta(E)-\delta(E_1)-\delta(E_2) = - {\alpha^2 \over 2rr_1r_2} 
\]
where $\alpha = r_1 \alpha_2 - r_2 \alpha_1$. 
Then Lemma \ref{bogomolovlemma} implies that 
\[
- {\alpha^2 \over 2rr_1r_2}\leq \delta(E).
\]
For simplicity let $\alpha = a \Xi + b f$, $a,b\in \IZ$.
The slope equality in \eqref{eq:extconditions} implies that 
$\alpha \cdot \omega|_{S'}=0$, which yields 
\[
b = 2a\left(1-{s\over t}\right).
\]
Therefore 
\[
-\alpha^2 = 2a^2\left({2s\over t}-1\right).
\]
Next note that $a\neq 0$; if $a=0$, one has $b=0$ as well, hence $\alpha=0$, 
contradicting the second condition in \eqref{eq:extconditions}. 
Therefore $a^2\geq 1$ since $a\in \IZ$. Moreover, 
$\delta(E)\geq 0$ according to Lemma \ref{bogomolovlemma}, and 
$1\leq r_1,r_2 \leq r$. This implies inequality \eqref{eq:wallboundA}. 

\hfill $\Box$

\begin{lemm}\label{midpointlemma} 
Let $E$ be an $\omega$-slope stable pure dimension two sheaf 
on $X$ with numerical invariants 
\[
\ch_1(E) = rD, \qquad \ch_2(E) = lf, \qquad \ch_3(E) = -n \ch_3(\CO_x),
\]
$l,n,r\in \IZ$, $r\geq 1$. 
Suppose there exists $t' \in \IR$, $0<t'<t$ such that $E$ is not 
$\omega'$-slope semistable, where $\omega' = t'\Theta -sp^*K_B$. 
Then 
\be\label{eq:wallboundB}
{t\over s} > {2\over 1+ r^3\delta(E)}.
\ee
\end{lemm} 

{\it Proof.} Any sheaf $E$ with  $\ch_1(E) =rD$ 
must be set theoretically supported on a finite union of 
$K3$ fibers of $X$. Since $E$ is $\omega$-slope stable, it must be 
scheme theoretically supported on a reduced irreducible 
fiber $\iota:S\hookrightarrow X$. 

Let $\CQ_E(t',t)$ denote the family 
of sheaves $E'$ such that $E'$ is a nonzero pure dimension 
two quotient of $E$, not isomorphic to $E$, and 
$\mu_{\omega'}(E') < \mu_{\omega'}(E)$. According to Grothendieck's 
lemma \cite[1.7.9]{huylehn}, $\CQ_E(t',t)$ is bounded. Any quotient $E'$ of $E$ is also scheme theoretically supported on $S$ and has invariants 
of the form 
\be\label{eq:quotinv}
\ch_1(E') = r'D, \qquad \ch_2(E') = m'\Xi+l'f, \qquad \ch_3(E') = -n' \ch_3(\CO_x),
\ee
$l',m',n',r'\in \IZ$, $r'\geq 1$. 
Since the family $\CQ_{E}(t',t)$ is bounded, the set of 
numerical invariants $(r',m',l',n')$ of all sheaves in this family is 
finite. 

For any $t''\in \IR$, $t'\leq t''\leq t$ set $\omega'' = t''\Theta -s p^*K_B$. 
For any $\gamma'=(l',m',n',r')\in \IZ^4$, $r'\geq 1$ let $\eta_{\gamma'}: [t',\ t] \to \IR$ be the linear function 
\[ 
\eta_{\gamma'}(t'') = {2m'\over r'} s - \left({2m'-l'\over r'} + {l\over r}\right)t''. 
\]
Then note that for any sheaf $E'$ with invariants \eqref{eq:quotinv} one has 
\[
\mu_{\omega''}(E') - \mu_{\omega''}(E) = 
{\eta_{\gamma'}(t'')\over t''(2s- t'')}.  
\]
Since $E$ is $\omega$-slope stable and not $\omega'$-slope semistable, one has 
\[
\eta_{\gamma'}(t') <0, \qquad \eta_{\gamma'}(t)>0 
\]
for any sheaf $E'$ in the family $\CQ_E(t',t)$. Therefore $\eta_{\gamma'}$ is an increasing linear function of $t''$ for any such sheaf . 
In particular there exists exactly one point $t'<t(\gamma') <t$ such that $\eta_{\gamma'}(t(\gamma'))=0$. 
The set of all $t(\gamma')$ associated to $E'$ in $\CQ_E(t',t)$ is finite. 
Let $t_0$ be its maximal element and $\omega_0 =t_0\Theta - sp^*K_B$.  Then it will be shown below 
that $E$ is strictly $\omega_0$-slope semistable. 

Given the choice of $t_0$, one 
has $\eta_{\gamma'}(t_0) \geq 0$ for any quotient $E\twoheadrightarrow
E'$ in $\CQ_E(t',t)$. Moreover, there exists $E_0'$ in $\CQ_E(t',t)$ 
such that $\eta_{\gamma'}(t_0)=0$. Clearly, $E_0'$ cannot be isomorphic 
to $E$ since $\mu_{\omega'}(E_0)<\mu_{\omega'}(E)$. Hence the kernel $E_0'' = 
{\rm Ker}(E\twoheadrightarrow E_0')$ is nontrivial. 
This implies that $\ch_1(E_0')=r'_0D$, $\ch_1(E_0'') = r_0''D$ 
with $r_0',r_0''\geq 1$. 

Given a quotient $E\twoheadrightarrow E'$ not in $\CQ_E(t',t)$, one has 
\[
\eta_{\gamma'}(t')\geq 0, \qquad \eta_{\gamma'}(t)>0.
\]
Since $\eta_{\gamma'}$ is linear this implies that $\eta_{\gamma'}(t_0) >0$, hence $E'$ cannot destabilize $E$ with respect to $\omega_0$. 

In conclusion $E$ is indeed $\omega_0$-slope semistable and 
there is an exact sequence 
\[
0\to E_0'' \to E\to E_0'\to 0 
\]
such that $\mu_{\omega_0}(E_0'') = \mu_{\omega_0}(E_0')$
and $r_0',r_0''\geq 1$. 
Moreover, since $E$ is $\omega$-slope stable one must have 
\[ 
{1\over r'_0} \ch_2(E'_0) - {1\over r_0''}\ch_2(E_0'') \neq 0.
\]
Then Lemma \ref{wallbound} implies that $t_0/s\geq {2/( 1+ r^3\delta(E))}$. 

\hfill $\Box$

\subsection{Proof of Proposition \ref{Kthreepencil} }

Let $(n,r)\in \IZ\times \IZ$ be fixed integers, $n\geq 0$, $r\geq 1$. 
For any $j\in \IZ$, $1\leq j\leq r$, let 
\[
\Gamma_j(n,r)= \big\{((n_1,r_1), \ldots (n_j,r_j)) \in (\IZ\times \IZ)^{\times j}\, |\, n_i\geq 0, \ r_i\geq 1,\ 1\leq i \leq j, \
\sum_{i=1}^j r_i =r, \ \sum_{i=1}^j n_i =n \big\}.
\]
Then let $\Gamma(n,r)= \cup_{j=1}^r \Gamma_j(n,r)$. 
Clearly $\Gamma(n,r)$ is a finite set. Let $t\in \IR$, $t>0$ be such 
that $t/s\in \IR\setminus \IQ$ and 
\be\label{eq:smallt}
{t\over s} <  {2\over 1+ r_i^3 n_i}, \qquad 1\leq i \leq j, 
\ee
for any element $\big((n_i,r_i)\big)_{1\leq i \leq j}\in \Gamma_j(n,r)$, 
and for all $1\leq j \leq r$. 

Let $E$ be an
$\omega$-semistable sheaf on $X$ with 
numerical invariants 
\[
\ch_1(E) = rD, \qquad \ch_2(E) = lf, \qquad \ch_3(E) = -n \ch_3(\CO_x),
\]
$l,n,r\in \IZ$, $r\geq 1$. Then $E$ is $\omega$-slope semistable. Let 
\be\label{eq:EJHfiltr}
0= E_0 \subset E_1 \subset \cdots \subset E_j = E 
\ee
be a Jordan-H\"older filtration of $E$ with respect to $\omega$-slope 
stability.  Let 
\be\label{eq:JHfactors} 
\ch_1(E_i/E_{i-1}) = r_i D, \qquad 
\ch_2(E_i/E_{i-1}) = m_i \Xi+l_i f, \qquad 
\ch_3(E_i/E_{i-1}) = -n_i \ch_3(\CO_x)
\ee 
be the numerical invariants of the $i$-th successive quotient, 
 where $r_i, l_i, n_i \in \IZ$, $r_i\geq 1$. 
Since $t/s \in \IR\setminus\IQ$ a simple computation 
shows that 
\be\label{eq:JHfactorsB} 
m_i =0, \qquad {l_i\over r_i}={l\over r} 
\ee 
for each $1\leq i\leq j$. 
Obviously, 
\[ 
\sum_{i=1}^j r_i =r, \qquad \sum_{i=1}^j n_i =n
\]
Moreover, Lemma \ref{bogomolovlemma} shows that $\delta(E_i/E_{i-1}) 
= n_i \geq 0$ for each $1\leq i\leq j$.
Since $t/s$ satisfies inequalities \eqref{eq:smallt}, Lemma \ref{wallbound} 
implies that each $E_i/E_{i-1}$ is adiabatically $\omega$-slope semistable. 
According to Lemma \ref{genstablemma}, this implies that each 
$E_i/E_{i-1}$ is generically semistable as in Definition \ref{genstabdef}. 
Let $H$ be a very ample divisor in $B$ satisfying the genericity conditions 
in loc. cit. for $E$ as well as for each successive quotient 
$E_i/E_{i-1}$. In particular $Z=p^{-1}(H)$ is a smooth elliptic surface which 
intersects the set theoretic support of $E$ along a finite union of elliptic 
fibers. Then 
Lemma \ref{zerotor} implies that the filtration \eqref{eq:EJHfiltr} 
restricts to a filtration of $E|_Z$ with successive quotients 
$(E_i/E_{i-1})|_Z$,  $1\leq i \leq j$, 
and $\chi((E_i/E_{i-1}))|_Z=0$ for all $1\leq i \leq j$. 
Since each $E_i/E_{i-1}$ is generically semistable, $(E_i/E_{i-1})|_Z$ is a zero slope 
semistable pure dimension one sheaf on $Z$. Hence $E|_Z$ is also 
semistable, which means that $E$ is generically semistable. Finally, Lemma 
\ref{adiabstablemma} implies that $E$ is adiabatically 
$\omega$-semistable. 

\hfill $\Box$

\bibliography{FM_ref.bib}
\bibliographystyle{abbrv}
\end{document}